\documentclass{amsart}
\usepackage{amsmath,amsfonts,epsfig}

\begin{document}

\newtheorem{thm}{Theorem}[section]
\newtheorem{lem}[thm]{Lemma}
\newtheorem{cor}[thm]{Corollary}
\newtheorem{mainlem}[thm]{Main Lemma}
\newtheorem{prop}[thm]{Proposition}
\newtheorem{conj}[thm]{Conjecture}

\theoremstyle{definition}
\newtheorem{defn}{Definition}[section]

\newcommand{\U}{\ensuremath{\widetilde}}
\newenvironment{pfsketch}{{\it Sketch of Proof:}\quad}{\square \vskip 12pt}

\newtheorem{rmk}[thm]{Remark}

\newcommand{\Hn}{\ensuremath{\mathbb{H}^n}}
\newcommand{\h}{\ensuremath{{\text{hyp}}}}
\newcommand{\g}{\ensuremath{{\text{geod}}}}
\newcommand{\Hess}{\ensuremath{{\text{Hess} \ }}}

\def\square{\hfill${\vcenter{\vbox{\hrule height.4pt \hbox{\vrule width.4pt
height7pt \kern7pt \vrule width.4pt} \hrule height.4pt}}}$}

\newenvironment{pf}{{\it Proof:}\quad}{\square \vskip 12pt}

\title{The barycenter method on singular spaces}
\author{Peter A. Storm}
\thanks{This research was partially supported an NSF Postdoctoral Fellowship.}
\date{June 7th, 2004}

\begin{abstract}
Compact convex cores with totally geodesic boundary are proven to uniquely minimize volume over all hyperbolic
$3$-manifolds in the same homotopy class.  This solves a conjecture in Kleinian groups concerning acylindrical
$3$-manifolds.  Closed hyperbolic manifolds are proven to uniquely minimize volume over all compact hyperbolic
cone-manifolds in the same homotopy class with cone angles $\le 2 \pi$.  Closed hyperbolic manifolds are proven to
minimize volume over all compact Alexandrov spaces with curvature bounded below by $-1$ in the same homotopy
class.  A version of the Besson-Courtois-Gallot theorem is proven for $n$-manifolds with boundary.  The proofs
extend the techniques of Besson-Courtois-Gallot.
\end{abstract}

\maketitle

\section{Introduction}\label{introduction}

This paper extends the barycenter map machinery of Besson-Courtois-Gallot \cite{BCGlong} to a class of singular
metric spaces called \emph{convex Riemannian amalgams} (defined in Section \ref{amalgams}).  This class of
singular spaces includes cone-manifolds and the metric doubling of hyperbolic convex cores across their boundary.
These singular space techniques are used to solve a conjecture in Kleinian groups due to Bonahon.  Specifically,
we prove that compact convex cores with totally geodesic boundary uniquely miminize volume over all hyperbolic
$3$-manifolds in the same homotopy class (see Theorem \ref{cores thm}).

The main tool of this paper is the following extension of Besson-Courtois-Gallot techniques to convex Riemannian
amalgams. \vskip 6pt \noindent\textbf{Theorem \ref{mainthm}.}  \itshape  For $n \ge 3$, let $Z$ be a compact
$n$-dimensional convex Riemannian amalgam.  Let $M_\h$ be a closed hyperbolic $n$-manifold.  Let $h(\U{Z})$ denote
the volume growth entropy of the universal cover of $Z$.  If $f: Z \longrightarrow M_\h$ is a homotopy equivalence
then $$h(\U{Z})^n \, \text{Vol}(Z) \ge (n-1)^n \text{Vol}(M_\h)$$ with equality if and only if $f$ is homotopic to
a homothetic homeomorphism. \normalfont

\vskip 6pt {\noindent}By restricting attention to cone-manifolds, we obtain

\vskip 6pt \noindent\textbf{Theorem \ref{cone-manifolds theorem}.}  \itshape  For $n \ge 3$, let $Z$ be compact
$n$-dimensional cone-manifold built with simplices of constant curvature $K \ge -1$.  Assume all its cone angles
are $\le 2\pi$.  Let $M_\h$ be a closed hyperbolic $n$-manifold.  If $f: Z \longrightarrow M_\h$ is a homotopy
equivalence then $$\text{Vol}(Z) \ge \text{Vol}(M_\h)$$ with equality if and only if $f$ is homotopic to an
isometry. \normalfont \vskip 6pt {\noindent}If $Z$ is allowed to be any Alexandrov space with curvature bounded
below by $-1$ (see Section \ref{Alexandrov section}), we then obtain

\vskip 6pt \noindent\textbf{Theorem \ref{Alexandrov theorem}.} \itshape Let $Z$ be a compact $n$-dimensional $(n
\ge 3)$ Alexandrov space with curvature bounded below by $-1$.   Let $M_\h$ be a closed hyperbolic $n$-manifold.
If $f: Z \longrightarrow M_\h$ is a homotopy equivalence then $$\text{Vol}(Z) \ge \text{Vol} (M_\h).$$ \normalfont
\vskip 6pt

As mentioned above, the paper's main theorems solve a conjecture in Kleinian groups.  To state things precisely,
let H$(N)$ denote the set of marked oriented isometry classes of hyperbolic $3$-manifolds $M$ equipped with a
homotopy equivalence $N \longrightarrow M$.  Define a volume function $$\text{Vol}: M \in \text{H}(N) \longmapsto
\text{Vol}(C_M).$$ It is a consequence of Thurston's Geometrization Theorem and Mostow Rigidity that $N$ is
acylindrical if and only if there exists a convex cocompact $M_g \in \text{H}(N)$ such that $\partial C_{M_g}$ is
totally geodesic \cite[pg.14]{Th2}.  Moreover, $M_g$ is unique up to isometry.  Let $M_g^{\text{opp}}$ denote
$M_g$ with the opposite orientation. \vskip 3pt {\noindent}\textbf{Conjecture:} $M_g$ and $M_g^{\text{opp}}$ are
the \emph{only} global minima of Vol over H$(N)$. \vskip 3pt {\noindent}Initial progress on this conjecture was
made by Bonahon \cite{Bon}.  Using different methods, Bonahon proved that $M_g$ is a strict local minimum of Vol
(in the quasi-isometric topology on H$(N)$).  In \cite{S}, the author proved that $M_g$ and $M_g^{\text{opp}}$ are
global minima of Vol.  Here this conjecture is completely solved by \vskip 6pt \noindent\textbf{Theorem \ref{cores
thm}.} \itshape Let $N$ be a compact acylindrical $3$-manifold.  Let $M_g \in \text{H}(N)$ be a convex cocompact
hyperbolic $3$-manifold such that the boundary of the convex core $\partial C_{M_g} \subset M_g$ is totally
geodesic.  Then for all $M \in \text{H}(N)$, $$\text{Vol}(C_M) \ge \text{Vol} (C_{M_g}),$$ with equality if and
only if $M$ and $M_g$ are isometric. \normalfont \vskip 6pt

The techniques used to prove Theorem \ref{cores thm} immediately generalize to a version of the
Besson-Courtois-Gallot theorem for manifolds with boundary.  (Perelman's unpublished Doubling theorem
\cite[Thm.5.2]{P} is used in the proof of Theorem \ref{n-mfds}.  See Theorem \ref{doubling theorem}.)

\vskip 6pt \noindent\textbf{Theorem \ref{n-mfds}.} \itshape Let $Z$ be a compact convex Riemannian $n$-manifold
with boundary $(n \ge 3)$.  Assume the sectional curvature of int$(Z)$ is bounded below by $-1$.  Let $Y_\g$ be a
compact convex hyperbolic $n$-manifold with totally geodesic boundary.   Let $f: (Z , \partial Z) \longrightarrow
(Y_\g, \partial Y_\g )$ be a homotopy equivalence of pairs.  Then $$\text{Vol} (Z) \ge \text{Vol} (Y_\g),$$ with
equality if and only if $f$ is homotopic to an isometry. \normalfont \vskip 6pt

This paper's method of proof may be of independent interest.  Following \cite{BCGlong}, these results are proven
by defining a \emph{natural map} from a nice path metric space $Z$ to a hyperbolic manifold $M_\h$.  Instead of
obtaining this map as a uniform limit of approximating maps, the natural map is here obtained in a single step.
The idea is to emulate the ``short'' proof of the Besson-Courtois-Gallot theorem found in \cite{BCGergodic}, where
it is additionally assumed that $Z$ is nonpositively curved.  Here the assumption of nonpositive curvature is
removed, but the gist of the ``short'' proof is retained.  As a cost for this generalization, the arguments here
require that $Z$ and $M_\h$ be homotopy equivalent.  In \cite{BCGlong}, only a map $Z \longrightarrow M_\h$ of
nonzero degree is required.

\vskip 6pt The results in this paper represent a large portion of the author's Ph.D. thesis, completed at the University of Michigan.  The author thanks his advisor, Richard Canary, for his essential assistance at every stage of this
project.  The author also enjoyed several helpful conversations with Yair Minsky.  Thanks to Ralf Spatzier for
introducing the author to the work of Besson-Courtois-Gallot.

\subsection{Sketch of the proof of the Main Theorem}

Let $Z$ be a compact $n$-dimensional convex Riemannian amalgam with universal cover $X$ (e.g. a cone-manifold or
the double of a convex core), $M_\h$ a closed hyperbolic $n$-manifold, and $f: Z \longrightarrow M_\h$ a homotopy
equivalence.  Up to rescaling the metric of $Z$, way may assume that $h(X)= (n-1) = h(\Hn)$.  The first goal is to
find a volume decreasing map $F: Z \longrightarrow M_\h$ homotopic to $f$.  The second goal is to show that if the
volume decreasing map $F$ is in fact volume preserving, then it is an isometry.

\vskip 3pt {\noindent}\textbf{Step 1:} \textit{Defining ``visual measures'' on $\partial X$. (Section
\ref{horoboundary and densities}).}  Following \cite{BM}, we define a generalization of Patterson-Sullivan measure
$\{\mu_x \}_{x \in X}$ supported on a function theoretic compactification $\mathcal{H} X$ of $X$ by Busemann
functions.  As $Z$ and $M_\h$ are homotopy equivalent, $X$ must be Gromov hyperbolic with Gromov boundary
$\partial X$.  We define an Isom$(X)$-equivariant continuous surjection $\pi: \mathcal{H}X \longrightarrow
\partial X$, and use it to push forward the measures $\{ \mu_x \}$ onto $\partial X$.  The resulting family of
probability measures $\{ \pi_* \mu_x \}$ are the ``visual measures'' used to define the natural map $F$.

\vskip 3pt {\noindent}\textbf{Step 2:} \textit{Defining $F$ (Section \ref{inequality}).} The homotopy equivalence
$f$ lifts to a quasi-isometry between the universal covers $f:X \longrightarrow \Hn$, which in turn induces a
homeomorphism of Gromov boundaries $f: \partial X \longrightarrow \partial \Hn$.  The measure $\pi_* \mu_x$ is
pushed forward via the boundary homeomorphism to obtain a measure $f_* \pi_* \mu_x$ on $\partial \Hn$.  Finally,
$F(x) \in \Hn$ is defined to be the barycenter of the measure $f_* \pi_* \mu_x$.  In sum,
\begin{eqnarray*}
F :  X      & \longrightarrow & \Hn  \\ x & \longmapsto &  \text{barycenter of measure } (f \circ \pi)_* \mu_x
\end{eqnarray*}
It is relatively easy to show that $F$ descends to a map $F: Z \longrightarrow M_\h$, and is homotopic to $f: Z
\longrightarrow M_\h$.  In particular, $F: Z \longrightarrow M_\h$ is surjective.

\vskip 3pt {\noindent}\textbf{Step 3:} \textit{$F$ is locally Lipschitz (Section \ref{lipschitz}).}  In order to
use calculus to study the map $F$, we must prove it is locally Lipschitz.  (For example, a locally Lipschitz map
which is infinitesimally volume decreasing almost everywhere must be volume decreasing.)  This is done by
factoring $F$ as a composition $F = P \circ \Phi$, such that $P$ and $\Phi$ can be analyzed directly.  $\Phi : X
\longrightarrow L^2_+ (\mathcal{H} X) \subset L^2 (\mathcal{H} X)$ equivariantly maps $X$ into the strictly
positive functions of the Hilbert space $L^2 (\mathcal{H} X)$.  $P : L^2_+ (\mathcal{H}X) \longrightarrow \Hn$ is
basically the barycenter map, thinking of $L^2_+ (\mathcal{H}X)$ as a space of measures.  We show $\Phi$ is
locally Lipschitz by direct estimates.  Applying the implicit function theorem shows $P$ is $\mathcal{C}^1$.
Together this shows $F = P \circ \Phi$ is locally Lipschitz.

\vskip 3pt {\noindent}\textbf{Step 4:} \textit{$F$ is infinitesimally volume decreasing a.e. (Section
\ref{jacobian}).}  With only minor modifications, the arguments of \cite{BCGergodic} can be applied to show $|
\text{Jac} F| \le 1$ almost everywhere.  This accomplishes the first goal of showing $F$ is volume decreasing.
The arguments of \cite{BCGergodic} also show that if $| \text{Jac} F (x)| =1$ for some $x$, then $dF_x$ is an
infinitesimal isometry.  Thus if $F$ is volume preserving, it must be an infinitesimal isometry almost everywhere.

\vskip 3pt {\noindent}\textbf{Step 5:} \textit{$F$ is volume preserving implies it is an isometry (Section
\ref{equality}).}  Applying the arguments from \cite[pg.790-793]{BCGlong}, we show a volume preserving map $F$ is
a local isometry on an open dense set.  We show $F$ is injective by using some local properties of convex
Riemannian amalgams.  Thus $F$ is a homeomorphism.  Again using the convex Riemannian structure on $Z$, we prove
$F$ is an isometry.  This accomplishes the second goal, and completes the proof of Theorem \ref{mainthm}.

\subsection{Sketch of the applications}
With the above machinery established, the theorems concerning cone-manifolds (Section \ref{cone-manifolds
section}), Alexandrov spaces (Section \ref{alexandrov spaces section}), and manifolds with boundary (Section
\ref{section n-mfds}) are easy to prove.  To apply the machinery to hyperbolic convex cores requires more work
(Section \ref{kleinian section}).

Recall the hypotheses of Theorem \ref{cores thm}.  Let $M$ and $M_g$ be homotopy equivalent acylindrical convex
cocompact hyperbolic $3$-manifolds.  Assume the convex core $C_{M_g} \subset M_g$ has totally geodesic boundary.
The goal is to prove Vol$(C_M) \ge \text{Vol} (C_{M_g})$, with equality if and only if $M$ and $M_g$ are
isometric.  To begin, metrically double the convex cores across their boundaries to obtain the convex Riemannian
amalgams $DC_M$ and $DC_{M_g}$.  Notice that $DC_{M_g}$ is in fact a closed hyperbolic $3$-manifold.  A short
argument proves there exists a $\pi_1$-equivariant quasi-isometric homeomorphism $f: \U{DC_M} \longrightarrow
\U{DC_{M_g}}$.  Applying the above machinery to $f$ yields Theorem \ref{cores thm}.

In the case where $M$ is geometrically finite with at least one rank one cusp, geometric (but not hyperbolic) Dehn
surgery arguments are used to reduce to the case of closed (non-hyperbolic) manifolds, where Theorem \ref{mainthm}
can be applied.  (These geometric Dehn surgery techniques are based on \cite{Bes,L}.)

%Extending this argument to include parabolics requires some additional notation and one standard fact.  Namely, performing hyperbolic Dehn surgery on a finite volume hyperbolic $3$-manifold always strictly decreases volume.

\section{Preliminaries}
The following is a review of the necessary definitions.  Throughout this paper, metric spaces are assumed to be
complete unless otherwise stated.

\subsection{$\delta$-Hyperbolicity}
This paper will follow the definitions and notation of \cite{GH}.  For convenience, we recall a few basic notions.
Let $(X,d)$ be a $\delta$-hyperbolic space with basepoint $o \in X$.  Then for $x,y \in X$, the \emph{Gromov
product} of $x$ and $y$ is $$( x \, | \, y) := \frac{1}{2} \left\{
    d(x,o) + d(y,o) - d(x,y) \right\}.$$
A defining property of $\delta$-hyperbolic spaces is that for any triple $x,y,z \in X$,
\begin{eqnarray}\label{3}
(x \, | \, y) &\ge & \text{min} \{(x \, | \, z), (z \, | \, y) \} - \delta.
\end{eqnarray}

\begin{figure}[ht]
\begin{center}
\epsfig{file=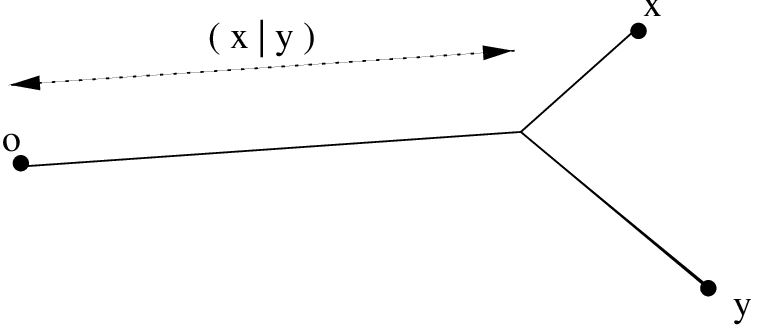,angle=-90} \label{gromov product} \caption{}
\end{center}
\end{figure}

{\noindent}The geometric content of the Gromov product may be difficult to grasp initially.  The idea is that
geodesic triangles in a $\delta$-hyperbolic space are very close to being tripods.  For a tripod the Gromov
product has the simple interpretation shown in Figure \ref{gromov product}: it is the length of the tripod's
``$o$'' leg.
{\noindent}Since long geodesic triangles in a $\delta$-hyperbolic space are very close to tripods, on a
sufficiently large scale the Gromov product is the length of the ``$o$'' leg of the geodesic triangle formed by
$o,x$, and $y$.

The Gromov boundary at infinity of $X$ will be denoted by $\partial X$.  The Gromov product can be extended to all
of $X \cup \partial X$ as follows.  For $a,b \in X \cup \partial X$, define $$( a \, | \, b ) := \sup \liminf_{i,j
\rightarrow \infty} ( x_i \, | \, y_j )$$ where the supremum is taken over all sequences $\{ x_i \}, \{y_j \}
\subset X$ such that $x_i \rightarrow a$ and $y_j \rightarrow b$.  For $a,b \in X$ this reduces to the previous
definition.  Using inequality (\ref{3}), one can show that for any given sequences $x_i \rightarrow a$, $y_j
\rightarrow b$, $$( a \, | \, b ) - 2 \delta \le \liminf (x_i \, | \, y_j) \le ( a \, | \, b)$$ (see
\cite[pg.122]{GH}).  Therefore a sequence $\{ z_i \} \subset X \cup \partial X$ converges to $a \in \partial X$ if
and only if $( z_i \, | \, a ) \longrightarrow \infty$.

\subsection{Barycenter}

Consider $\Hn$ with basepoint $o$.  For each $\theta \in \partial \Hn$ let $B^o_\theta$ be the unique Busemann
function of $\theta$ such that $B^o_\theta (o) = 0$.  Namely, if $\gamma$ is the geodesic ray based at $o$
asymptotic to $\theta$, then $B^o_\theta (p) = \lim_{t \rightarrow \infty} [d(p, \gamma(t)) - t]$.  Let $\lambda$
be a Radon measure on the compact space $\partial \Hn $.  Define the average $$B_\lambda (p):= \int_{\partial \Hn}
B^o_\theta (p) \ d\lambda (\theta).$$

\begin{prop}
\label{positive definite} \cite[Appendix A]{BCGlong} If $\lambda$ has no atoms, then $B_\lambda$ is proper and has
a unique critical point in $\Hn$ corresponding to the unique global minimum.  Moreover, the Hessian of $B_\lambda$
is a positive definite bilinear form.  Namely, for all $v \in T_p \Hn$, $$(\Hess B_\lambda) (v,v) :=
    \langle \nabla_v \nabla B_\lambda, \, v \rangle > 0.$$
\end{prop}
The unique critical point of $B_\lambda$ is the \textit{barycenter} of $\lambda$, denoted $\text{bar}\, \lambda$.
Since the barycenter is the unique critical point of $B_\lambda$, it is defined implicitly as the unique point $p$
such that $$\int_{\partial \Hn} \langle \nabla B^o_\theta , \ v
    \rangle_p  \, d\lambda (\theta) = 0$$
for all $v \in T_p \Hn$.  Notice that the barycenter map is scale invariant, i.e. the barycenter of the measure
$\lambda$ equals the barycenter of the measure $c \lambda$ for any $c>0$.

\subsection{The brain in a jar lemma}
At a crucial stage in the proof of Proposition \ref{mainprop}, the following linear algebra lemma is needed.

\begin{lem}
\label{brain in a jar lemma} Let $H$ be an $n \times n$ positive definite symmetric matrix with trace$(H)=1$.  If
$n \ge 3$, then $$\frac{\det(H)}{[ \det (\text{Id} - H) ]^2} \le
    \left[ \frac{n}{(n-1)^2} \right]^n.$$
Moreover, equality holds if and only if $H = \frac{1}{n} \text{Id}$.
\end{lem}
{\noindent}This lemma is false for $n=2$.  It is the only part of this paper (and Besson-Courtois-Gallot theory in
general) which fails in $2$ dimensions.

\subsection{Generalized differentiable and Riemannian structures \normalfont{\cite{OS}}} \label{gdrs}
Let $X$ be a topological space, $\Omega \subseteq X$, and $n \in \mathbb{N}$.  A family $\{ (U_\phi, \phi
)\}_{\phi \in \Phi}$ is called a \emph{$\mathcal{C}^1$-atlas on $\Omega \subseteq X$} if the following hold:\\ (1)
For each $\phi \in \Phi$, $U_\phi$ is an open subset of $X$.\\ (2)  Each $\phi \in \Phi$ is a homeomorphism from
$U_\phi$ into an open subset of $\mathbb{R}^n$.\\ (3)  $\{ U_\phi \}_{\phi \in \Phi}$ is a covering of $\Omega$.\\
(4)  If two maps $\phi, \psi \in \Phi$ satisfy $U_\phi \bigcap U_\psi \not= \emptyset$, then $$\psi \circ
\phi^{-1} : \phi(U_\phi \bigcap U_\psi) \longrightarrow
        \psi (U_\phi \bigcap U_\psi)$$

is $\mathcal{C}^1$ on $\phi(U_\phi \bigcap U_\psi \bigcap \Omega)$.
\\

A family $\{g_\phi \}_{\phi \in \Phi}$ is called a \emph{$\mathcal{C}^0$-Riemannian metric} associated with a
$\mathcal{C}^1$-atlas $\{ (U_\phi, \phi) \}_{\phi \in \Phi}$ on $\Omega \subseteq X$ if the following hold:\\ (1)
For each $\phi \in \Phi$, $g_\phi$ is a map from $U_\phi$ to the set of positive symmetric matrices.\\ (2)  For
each $\phi \in \Phi$, $g_\phi \circ \phi^{-1}$ is continuous on $\phi(U_\phi \bigcap \Omega)$.\\ (3)  For any $x
\in U_\phi \bigcap U_\psi, \phi,\psi \in \Phi$, we have $$g_\psi (x) = [ d(\phi \circ \psi^{-1})(\psi(x))]^t \,
g_\phi (x) \,
        [d (\phi \circ \psi^{-1}) (\psi(x))].$$

These two structures induce a distance metric $D_g$ on $\Omega$ which we now describe.  The length of a piecewise
$\mathcal{C}^1$ path $\gamma: (0,1) \longrightarrow \Omega$ is defined in the usual way by pulling back the metric
tensor via $\gamma$.  A path $\eta : [0,1] \longrightarrow X$ is \emph{admissible} if $\eta^{-1} (X \setminus
\Omega)$ is a finite set of points $\{ t_1, \ldots ,t_l \}$, and $\eta$ is piecewise $\mathcal{C}^1$ on $(0,1)
\setminus \{ t_1, \ldots, t_l \}$.  For $x,y \in \Omega$, define $$D_g (x,y) := \inf \{ \text{ length} (\gamma) \
| \ \gamma \text{ joins } x \text{ to } y
    \text{ and is admissible} \}.$$
If $x$ and $y$ cannot be joined by an admissible path, set $D_g (x,y) = \infty$.  In general, the topology of
$(\Omega, D_g)$ can be quite different from the subspace topology of $\Omega \subseteq X$.

\subsection{Almost everywhere Riemannian spaces} \label{aers}
Let $(X,d)$ be a geodesic metric space with Hausdorff dimension $n < \infty$.  $X$ is an \emph{almost everywhere
Riemannian} metric space if there exists $\Omega \subseteq X$, a dense subset of full $n$-dimensional Hausdorff
measure, such that: \\ (1) $\Omega$ admits a $n$-dimensional $\mathcal{C}^1$-atlas $\{ (U_\phi, \phi )\}_{\phi \in
\Phi}$.  \\ (2) $\Omega$ admits a $\mathcal{C}^0$-Riemannian metric $\{ g_\phi \}_{\phi \in \Phi}$. \\ (3) Each
homeomorphism $\phi \in \Phi$ is in fact locally bilipschitz. \\ (4) The identity map $(\Omega, D_g)
\longrightarrow (\Omega,d)$ is an isometry (see Section \ref{gdrs}).\\ (5) The Riemannian metric induces a volume
element $d\text{vol}_\Omega$ on $\Omega$.  The measure on $X$ obtained by integrating this element equals
$n$-dimensional Hausdorff measure on $X$ ($X \setminus \Omega$ has zero measure).

Notice $X$ is not assumed to be a topological manifold, and $\Omega \subseteq X$ is not assumed to be open.
Conditions (1)-(5) are not as restrictive as they initially appear.  For example, Otsu and Shioya proved that a
finite dimensional Alexandrov space with curvature bounded below by $k \in \mathbb{R}$ is almost everywhere
Riemannian \cite{OS} (see Section \ref{Alexandrov section}).  Condition (3) allows Rademacher's theorem to be
applied to almost everywhere Riemannian spaces.  This means locally Lipschitz functions on an almost everywhere
Riemannian metric space are differentiable almost everywhere.  Condition (4) says that the
$\mathcal{C}^0$-Riemannian metric on $\Omega$ reproduces the metric $d$, in the sense that the metric completion
of $(\Omega, D_g)$ is isometric to $(X,d)$.

\subsection{Cone-manifolds \normalfont{\cite[pg.53]{CHK}}} \label{cone-manifolds}
An \emph{$n$-dimensional cone-manifold} $M$ is a manifold which can be triangulated so that the link of each
simplex is piecewise linear homeomorphic to a standard sphere and $M$ is equipped with a path metric such that the
restriction of the metric to each simplex is isometric to a geodesic simplex of constant curvature $K$.  The
singular locus $\Sigma$ consists of the points with no neighborhood isometric to a ball in a Riemannian manifold.

It follows that

{\noindent}$\bullet \ \  \Sigma$ is a union of totally geodesic closed simplices of dimension $n-2$.

{\noindent}$\bullet \ $   At each point of $\Sigma$ in an open $(n-2)$-simplex, there is a \emph{cone angle} which
is the sum of dihedral angles of $n$-simplices containing the point.

In particular, the singular locus of a $3$-dimensional cone-manifold forms a graph in the manifold.  A
cone-manifold is an almost everywhere Riemannian metric space (see Lemma \ref{cra is an aers}).  Though a
definition will not be given here, a cubed-manifold is another example of a cone-manifold (see \cite{AMR}).  More
abstractly, any manifold admitting a locally finite decomposition into convex geodesic polyhedra is a
cone-manifold.  This can be seen by adding superfluous faces to the polyhedral decomposition.

The manifold structure of a cone-manifold will not be used in this paper.  Theorem \ref{cone-manifolds theorem} is
equally valid for more general simplicial metric spaces not satisfying the above link condition.

\subsection{Convex Riemannian manifolds with boundary} \label{crmwb}
A geodesic metric space $C$ is an $n$-dimensional \emph{convex Riemannian} (resp. \emph{hyperbolic})
\emph{manifold with boundary} if \\ (1) $C$ is topologically an $n$-manifold with boundary,\\ (2) there is an
incomplete Riemannian (resp. hyperbolic) metric on the interior of $C$,\\ (3) the metric on $C$ is the metric
completion of the Riemannian (resp. hyperbolic) metric on the interior,\\ (4) for any pair of points in the
interior of $C$, the shortest path between them lies in the interior of $C$, and\\ (5) for any compact $K
\subseteq C$, the curvature of the Riemannian manifold $K \cap \text{int}(C)$ is bounded from above and below by
finite constants. \\ {\noindent}(Notice no differentiability assumptions have been made on the boundary.)

Property (5) ensures $C$ has ``locally bounded geometry''.  The lower curvature bound on compact sets guarantees
that $C$ is locally an Alexandrov space with curvature bounded below.  This implies the boundary $\partial C$ has
Hausdorff dimension $n-1$ \cite{OS}.  (In \cite{OS}, they prove the singular set of an $n$-dimensional Alexandrov
space with curvature bounded below has Hausdorff dimension $\le n-1$.  In this case, the boundary \emph{is} the
singular set.)  The local upper curvature bound gives a \emph{lower} bound on the volume of small metric balls in
$C$.  This useful property will be used in Section \ref{equality}.

Given a convex Riemannian manifold with boundary $C$, we can metrically double it across its boundary to obtain a
metric space $DC$.  Topologically $DC$ is the closed manifold obtained by doubling $C$ across its boundary.  The
metric on $DC$ is the path metric induced by gluing the two copies of $C$ (one with opposite orientation) along
$\partial C$.  The convexity of $C$ (property (4)) insures that the path metric obtained after doubling does not
alter the original metric on $C$.  Notice that $DC$ is an almost everywhere Riemannian metric space (see Lemma
\ref{cra is an aers}).

\subsection{Convex Riemannian amalgams} \label{amalgams}
A geodesic metric space $Z$ is an $n$-dimensional \emph{convex Riemannian amalgam} if $Z$ contains an
isometrically embedded locally finite countable collection $\{ C_j \} \subseteq Z$ of $n$-dimensional convex
Riemannian manifolds with boundary such that \\ (1) $\bigcup_j  C_j = Z, $ \\ (2)  int$(C_j) \cap \text{ int}(C_k)
= \emptyset$ for $j \neq k$.

\vskip 3pt {\noindent}(Notice $Z$ is not assumed to be a manifold.)  A cone-manifold is a convex Riemannian
amalgam (see Section \ref{cone-manifolds}).  Another convex Riemannian amalgam is the metric doubling $DC$ of a
convex Riemannian manifold with boundary $C$ (see Section \ref{crmwb}).

\begin{lem} \label{cra is an aers}
A convex Riemannian amalgam is an almost everywhere Riemannian metric space.
\end{lem}
\begin{pf}
Define $\Omega := \bigcup_j \text{int}(C_j)$.  We must check conditions (1)-(5) of Section \ref{aers}.  (1)-(3)
are trivial.  For (4) use the following consequence of convexity: any $x,y \in Z$ can be joined by a path $\gamma$
such that\\ $\bullet $  the length of $\gamma$ is arbitrarily close to $d(x,y)$, and \\ $\bullet $  $\gamma \cap
\partial C_j$ is at most two points for any $C_j$ in the decomposition.\\ For (5) use that $\bigcup_j \partial
C_j$ has measure zero.
\end{pf}
{\noindent}Convex Riemannian amalgams seem to be the most natural class of metric spaces for which the arguments
of Theorem \ref{mainthm2} are valid.

\subsection{Convex Cores} \label{convex cores}
Let $M$ be a hyperbolic manifold.  Let $S \subseteq M$ be the union of all closed geodesics in $M$.  The
\textit{convex core}, $C_M$, is the smallest closed convex subset of $M$ which contains $S$, in other words it is
the closed convex hull of $S$ in $M$.  The convex core may also be defined as the smallest closed convex subset of
$M$ such that the inclusion map is a homotopy equivalence.  $M$ is \emph{geometrically finite} if an
$\varepsilon$-neighborhood of $C_M$ has finite volume.  Otherwise, $M$ is \emph{geometrically infinite}.

For finite volume hyperbolic manifolds, the convex core is the entire manifold.  Thus this is a useful object only
in the infinite volume case, where $C_M$ is the smallest submanifold which carries all the geometry of $M$.

\subsection{Pared 3-manifolds \normalfont{\cite[Def.4.8]{M}}} \label{pared 3-manifolds section}
Let $N$ be a compact orientable irreducible $3$-manifold with nonempty boundary.  Assume $N$ is not a $3$-ball.
Let $P \subseteq \partial N$.  $(N,P)$ is a \emph{pared} $3$-manifold if the following three conditions hold. \\
(1)  Every component of $P$ is an incompressible torus or a compact annulus.\\ (2)  Every noncyclic abelian
subgroup of $\pi_1 (N)$ is conjugate into the fundamental group of a component of $P$.\\ (3)  Every
$\pi_1$-injective cylinder $C : (S^1 \times I, S^1 \times \partial I) \longrightarrow (N,P)$ is relatively
homotopic to a map $\psi$ such that $\psi (S^1 \times I) \subseteq P$.

By Thurston's Geometrization Theorem \cite{M}, $(N,P)$ is a pared $3$-manifold if and only if there exists a
geometrically finite hyperbolic structure on the interior of $N$ such that $C_M \cong N \setminus P$.

\subsection{Pared acylindrical \normalfont{\cite[pg.244]{Th3}}}
A pared $3$-manifold $(N,P)$ is \emph{pared acylindrical} if \mbox{$\partial N \setminus P$} is incompressible and
if every $\pi_1$-injective cylinder $$C: (S^1 \times I, S^1 \times \partial I) \longrightarrow (N,\partial N
\setminus P)$$ is homotopic rel boundary to $\partial N$.

\subsection{Deformation theory} \label{deformation theory section}
Let $\mu > 0$ be the Margulis constant for hyperbolic $3$-manifolds.  Then for a hyperbolic $3$-manifold $M$, the
$\mu$-thin part of $M$ is a disjoint union of bounded Margulis tubes and unbounded cusps \cite{BP}.  After
possibly making $\mu$ smaller, we may also assume that the intersection of $\partial C_M$ and the $\mu$-thin part
of $M$ is totally geodesic \cite[Lem.6.9]{M}.  Define $M^o$ to be $M$ minus the unbounded components of its
$\mu$-thin part.  In other words, $M^o$ is the manifold with boundary obtained by removing the cusps from $M$.

Let $(N,P)$ be a compact pared $3$-manifold.  Define the deformation space H$(N,P)$ as follows.  For a hyperbolic
$3$-manifold $M$ and a map $m: N \longrightarrow M^o$, $(M,m) \in \text{H}(N,P)$ if there exists a union $Q_M$ of
components of $\partial M^o$ such that $m: (N,P) \longrightarrow (M,Q_M)$ is a relative homotopy equivalence.
$(M_1, m_1) = (M_2, m_2)$ in H$(N,P)$ if there exists an orientation preserving isometry $\mathcal{I} : M_1
\longrightarrow M_2$ such that $\mathcal{I} \circ m_1 \sim m_2$.  (H$(N,P)$ admits several interesting topologies
\cite{Th3}.  They will not be needed here.)  $(M,m)$ has \emph{no additional parabolics} if $Q_M = \partial M^o$.
Using the product structure on the complement of $C_M$ \cite{EM} and the thick-thin decomposition, there exists a
relative homotopy equivalence $p:(M,Q_M) \longrightarrow (C_M \cap M^o, C_M \cap Q_M)$.

\begin{thm} \label{homeo thm} \emph{(Johannson)}
Let $(N,P)$ be a compact pared acylindrical $3$-manifold.  If $(M,m) \in \text{H}(N,P)$ is geometrically finite
then $p \circ m$ is homotopic to a homeomorphism $(N,P) \longrightarrow (C_M \cap
M^o, C_M \cap Q_M)$.
\end{thm}
\begin{pf}
Since $N$ is not homotopy equivalent to a surface, $C_M$ is a $3$-manifold.  Since $M$ is geometrically finite,
$(C_M \cap M^o, C_M \cap Q_M)$ is a compact pared $3$-manifold \cite[Cor.6.10]{M}.  So $p \circ m$ is a relative
homotopy equivalence between compact pared $3$-manifolds, and the domain is pared acylindrical.  By the work of
Johannson \cite[Lem.X.23,pg.235]{J}, $p \circ m$ is homotopic (rel the paring) to a homeomorphism $(N,P)
\longrightarrow (C_M \cap M^o, C_M \cap Q_M)$.
\end{pf}

We are interested in pared acylindrical $3$-manifolds because of the following corollary of Thurston's
Geometrization Theorem and Mostow rigidity.
\begin{cor} \label{pared acylindrical} \cite[pg.14]{Th2}
Let $(N,P)$ be a pared acylindrical $3$-manifold.  Then there exist exactly two spaces $(M,m), (M^{\text{opp}},
m^{\text{opp}}) \in H(N,P)$ such that $M$ and $M^{\text{opp}}$ are geometrically finite, $M$ and $M^{\text{opp}}$
have no additional parabolics, and the convex cores $C_M$ and $C_{M^{\text{opp}}}$ have totally geodesic boundary.
Moreover, there exists an orientation reversing isometry $\mathcal{I} : M \longrightarrow M^{\text{opp}}$ such
that \ $m^{\text{opp}} \sim \mathcal{I} \circ m$.
\end{cor}

\subsection{Alexandrov space with curvature bounded below by $-1$} \label{Alexandrov section}
There are many equivalent definitions of Alexandrov spaces.  Here we give the most common definition.  (See
\cite{BBI} for more information.)

Let $Y$ be a path metric space.  $Y$ is an \emph{Alexandrov space with curvature bounded below by $-1$} if about
each point in $Y$ there exists a neighborhood $U$ satisfying the following comparison condition.  Let $x,y,z \in
U$ be distinct points, let $w$ lie on the interior of a geodesic path $xy$ connecting $x$ to $y$.  Let $\tilde{x},
\tilde{y}, \tilde{z}, \tilde{w} \in \mathbb{H}^2$ be such that $d(x,y)=d(\tilde{x}, \tilde{y}), d(y,z)=
d(\tilde{y}, \tilde{z}), d(x,z) = d(\tilde{x}, \tilde{z}), d(x,w)=d(\tilde{x}, \tilde{w}), d(w,y) = d(\tilde{w},
\tilde{y})$.  We then require that $d(z,w) \ge d(\tilde{z}, \tilde{w})$.

\begin{figure}[ht]
\begin{center}
\epsfig{file=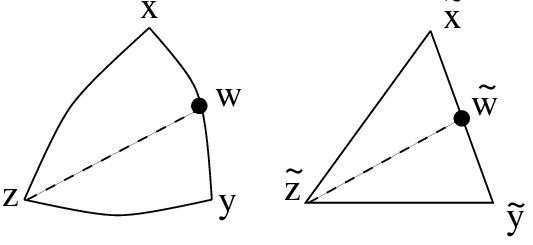,angle=-90} \caption{} \label{Alexandrov}
\end{center}
\end{figure}

This comparison condition guarantees that geodesic triangles in $Y$ are at least as fat as hyperbolic triangles.
The \emph{dimension} of an Alexandrov space with curvature bounded below is defined to be its Hausdorff dimension.
For finite dimensional Alexandrov spaces with curvature bounded below, the Hausdorff and topological dimensions
agree \cite[pg.21]{BGP}.

In Section \ref{alexandrov spaces section}, we will use
\begin{thm} \label{Otsu-Shioya} \cite{OS}
If $Y$ is a finite dimensional Alexandrov space with curvature bounded below by $-1$, then $Y$ is an almost
everywhere Riemannian metric space.
\end{thm}
\begin{pf}
Following \cite{OS}, let $S$ be the singular set of $Y$.  Define $\Omega := Y \setminus S$.  $\Omega$ is a
countable intersection of dense open sets with full measure.  It is therefore dense and has full measure.
Conditions (1)-(5) of Section \ref{aers} also follow from results in \cite{OS}.  Specifically, conditions (1) and
(2) follow from [Thm.B,pg.630], condition (3) follows from [Lem.5.1.3,pg.651], condition (4) follows from
[Thm.B,pg.630] and [Thm.6.4,pg.654], and condition (5) follows from [Sec.7.1].
\end{pf}

The main property of these Alexandrov spaces we will use is an upper bound on their volume growth entropy, which
we now define.

\subsection{Volume growth entropy}
Let $X$ be a geodesic metric space of Hausdorff dimension $n$, $\widetilde{X}$ be the universal cover of $X$, and
$\mathcal{H}^n$ be $n$-dimensional Hausdorff measure.  The \textit{volume growth entropy} of $\U{X}$ is the number
$$h(\U{X}) :=  \limsup_{R \rightarrow \infty}
   \frac{1}{R} \log \mathcal{H}^n ( B_{\widetilde{X}} (x,R)),$$
where $x$ is any point in $\widetilde{X}$, and the ball $B_{\widetilde{X}} (x,R)$ is in $\widetilde{X}$. \vskip
3pt {\noindent}(The volume growth entropy is independent of the choice of $x \in \widetilde{X}$.)  The following
theorem of Burago, Gromov, and Perelman will be important for this paper.

\begin{thm}
\label{Perel'man} \cite[pg.40]{BGP} If $X$ is an Alexandrov space with curvature bounded below by $-1$ and
Hausdorff dimension $n$, then the volume growth entropy of $\U{X}$ is less than or equal to the volume growth
entropy of $\mathbb{H}^n$.  In other words $$h(\U{X}) \le h(\mathbb{H}^n) = n-1.$$
\end{thm}

\section{Horoboundary and Densities}\label{horoboundary and densities}

In order to define a density in the necessary generality, a second definition of the boundary at infinity must be
made.  This is a more general definition than of the Gromov boundary at infinity; it makes sense for any proper
metric space.  In the case of negatively curved Riemannian manifolds, it reduces to the Gromov boundary.  (For
more information, see \cite[pg.21]{BGS}, \cite{BM}, and \cite[Sec.2]{F}.)

Let $Y$ be a proper metric space.  For $p \in Y$, define $d_p (y) := d(p,y)$.  Denote the space of continuous
(real valued) functions on $Y$ by $C(Y)$, and endow this set with the topology of uniform convergence on compact
sets.  Define an equivalence relation on $C(Y)$ by $f \sim g$ if and only if $f - g$ is a constant function.
Denote the quotient space $C(Y) / \sim$  by $C_* (Y)$.  $C_* (Y)$ is Hausdorff.  Define a map $\iota : Y
\longrightarrow C_* (Y)$ by $\iota (p) := [d_p]$.  $\iota$ is a topological embedding.

\begin{defn}
Let $\text{Cl} (Y)$ denote the closure of $\iota(Y)$ in $C_* (Y)$.  The \textit{horoboundary} of $Y$ is $$
\mathcal{H} Y := \text{Cl} (Y) \setminus \iota (Y).$$ A continuous function $h \in C(Y)$ such that $[h] \in
\mathcal{H} Y$ is a \textit{horofunction} of $Y$.

For $\eta \in C_* (Y)$ define a function $b_\eta : Y \times Y \longrightarrow \mathbb{R}$ by $$b_\eta (p,q) :=
h(p) - h(q) \text{ for any } h \in
     C(Y) \text{ such that } [h] = \eta.$$
If $Y$ has a fixed basepoint $o$, then define $b_\eta (p) := b_\eta (p,o)$.
\end{defn}

It is a quick check to see that $b_\eta$ is well-defined, i.e. independent of the choice of $h$.  The functions
$b_\eta$ are $1$-Lipschitz.  Thus applying the Arzela-Ascoli theorem shows Cl$(Y)$ is compact, implying
$\mathcal{H}Y$ is compact.  (If $Y$ is nonpositively curved, then horofunctions and the horoboundary are identical
to Busemann functions and the boundary at infinity \cite[pg.22]{BGS}.)

Isometries of $Y$ extend to homeomorphisms of $\mathcal{H} Y$ in the following simple manner.  Consider the
Isom$(Y)$-action by homeomorphisms on $C(Y)$ given by $\phi . f := f \circ \phi^{-1}$.  This action  descends to
an action on $C_* (Y)$.  Since $\phi . [d_p] = [d_{\phi.p}]$, the map $\phi |_{\iota (Y)} : \iota(Y)
\longrightarrow \iota (Y)$ is a homeomorphism.  Since $\phi$ acts as a homeomorphism on both $\iota (Y)$ and $C_*
(Y)$, it also acts as a homeomorphism on $\mathcal{H} Y$.  Thus we have defined an Isom$(Y)$-action by
homeomorphisms on $\mathcal{H} Y$.

\begin{defn}
Let $Y$ be proper metric space.  Let $G$ be a closed subgroup of $\text{Isom} (Y)$.  A continuous map (under the
weak-$*$ topology on measures) $$\mu : Y \longrightarrow \{ \text{positive Radon measures on } \mathcal{H} Y \}$$
is an $\ell$-\textit{dimensional density} for $G$ if \\ (1)  $\mu$ is $G$-equivariant, i.e. $\mu_{g . x} = g_*
\mu_x$, \\ (2)  $\mu_p \ll \mu_q$ for all $p,q$, and for $\eta \in \mathcal{H} Y$,
    $$\frac{d \mu_p}{d \mu_q} (\eta) = e^{- \ell (b_\eta (p,q))}.$$
\end{defn}

Example:  For $\Hn$, $\mathcal{H} \Hn = \partial \Hn $.  Let $\Phi: S_p \Hn \longrightarrow \partial \Hn$ be the
standard radial homeomorphism between the unit tangent sphere and the boundary at infinity.  Define $\mu_p$ to be
the push-forward of Lebesgue measure on $S_p Y$ by $\Phi$.  This is known as the \textit{visual measure at} $p$,
and is an $(n-1)$-dimensional density for all of Isom($\mathbb{H}^n$).  Visual measure is the most natural density
in this case.

\vskip 6pt Our entire reason for defining the horoboundary is the following proposition.

\begin{prop}
\label{PS measure} \cite[Prop. 1.1]{BM} Let $X$ be a proper metric space of Hausdorff dimension $n$ with basepoint
$o \in X$.  Let $m$ be $n$-dimensional Hausdorff measure.  Assume $m(X) = \infty$ and $X$ has finite volume growth
entropy $h(X)$.  Then there exists an $h(X)$-dimensional density $x \longmapsto \mu_x$ for Isom$(X)$.  This
density is called \textit{Patterson-Sullivan measure}.
\end{prop}
{\noindent}By normalizing we may always assume that $\mu_o$ is a probability measure.  Notice this normalization
implies that $\mu_{g.o} = g_* \mu_o$ is also a probability measure for any $g \in \text{Isom}(X)$.

\vskip 12pt Let $X$ be a $\delta$-hyperbolic proper path metric space.  Assume that $h(X) \in (0,\infty)$.  We
have defined two different compactifications of $X$; the Gromov boundary at infinity $\partial X$ and the
horoboundary $\mathcal{H} X$.  In general, those two compactifications are not homeomorphic.  However, they are
both necessary for the work of this paper.  To connect the two compactifications, we now define a continuous
Isom$(X)$-equivariant surjection $\pi : \mathcal{H} X \longrightarrow \partial X$.

Fix a basepoint $o \in X$.  Pick $\xi \in \mathcal{H} X$ and a sequence $\{ p_n \}$ such that $[ d_{p_n} ]
\rightarrow \xi$ in $C_* (X)$.  Then $\{ p_n \}$ leaves every compact set of $X$.  So by the Arzela-Ascoli theorem
there exists a subsequence $\{ a_l \} \subseteq \{ p_n \}$ such that the geodesic segments $oa_l$ converge to some
$[\gamma] \in \partial X$ in the compact-open topology, where $\gamma$ is a geodesic ray based at $o$.

Define $\pi: \mathcal{H}X \longrightarrow \partial X$ by
    $\pi(\xi) = [\gamma].$  Before proving $\pi$ is well defined, we need the following lemma.

\begin{lem}
In the above notation, the sequence $\{ a_l \}$ converges to $[ \gamma ] \in \partial X$.
\end{lem}

\begin{pf}
It is enough to show $(\, [\gamma] \, | \, a_l ) \longrightarrow \infty$.  Pick a large $M \gg 0$, and consider
the metric ball $B(o,M)$.  Find $L$ such that for $l \ge L$, the geodesic segment $oa_l$ is $1$-close to $\gamma$
on $B(o,M)$.  Then the picture looks roughly like Figure \ref{a}.

\begin{figure}[ht]
\begin{center}
\epsfig{file=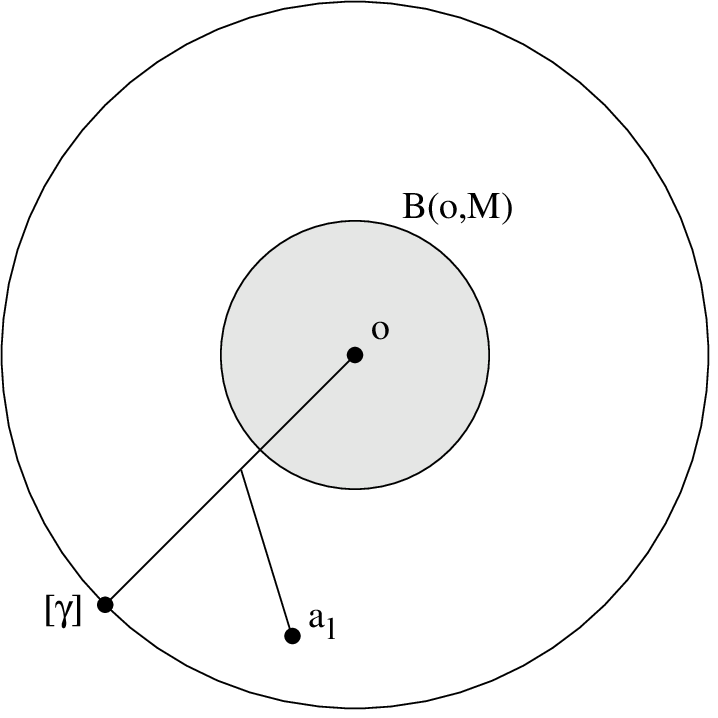,angle=-90} \caption{} \label{a}
\end{center}
\end{figure}

As $M \rightarrow \infty$, intuitively $( [\gamma] | a_l ) \longrightarrow \infty$.  To prove this carefully,
notice that for $k > 1$ $$d( \gamma ( kM), a_l) \le
    [d(a_l, o) - M] + [d( \gamma ( kM), o) - M ] + 1.$$
This implies $$( \gamma (kM) \, | \, a_l ) \ge M - 1/2.$$ So finally $$( [\gamma] \, | \, a_l ) \ge \liminf_{k
\rightarrow \infty}
    ( \gamma (kM) \, | \, a_l ) - 2 \delta
        \ge M - 2\delta - 1/2.$$
\end{pf}

\begin{lem}
The map $\pi: \mathcal{H} X \longrightarrow \partial X$ is well defined.
\end{lem}
\begin{pf}
Suppose $\pi$ is not well defined.  Then there exist sequences $ p_i, q_i  \longrightarrow \xi \in \mathcal{H}X$
such that $\{ p_i \}, \{q_i \}$ do not converge to a common point in $\partial X$.  Thus there exists an $M >0$
such that $( p_i \, | \, q_i ) < M$ for all $i$.

Pick a metric closed ball $K$ much larger than $B(o,M)$.  $ p_i, q_i  \longrightarrow \xi \in \mathcal{H}X$
implies that for large $i$ and all $x \in K$,
\begin{equation} \label{equation (*)}
d(p_i, x) - d(p_i, o) \approx d(q_i, x) - d(q_i , o).
\end{equation}
This will lead to a contradiction.

\begin{figure}[ht]
\begin{center}
\epsfig{file=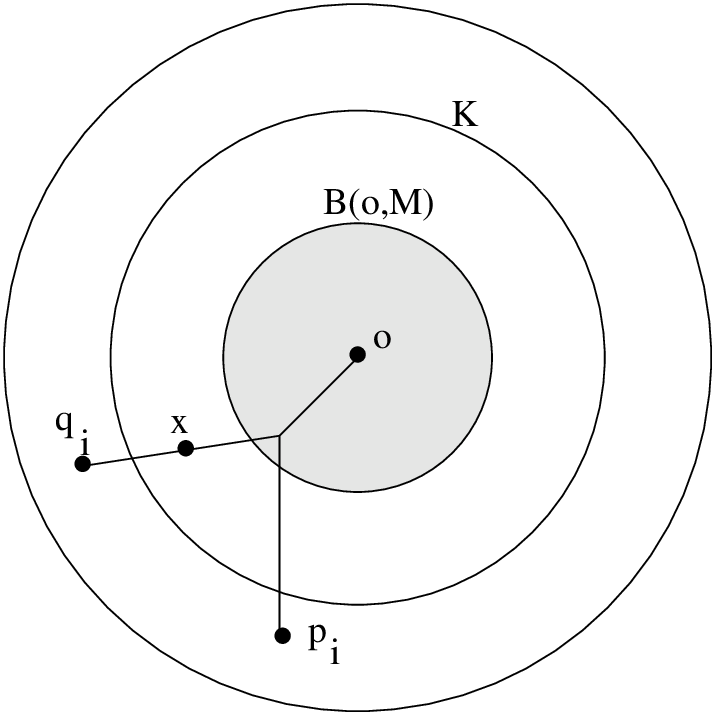,angle=-90} \caption{} \label{b}
\end{center}
\end{figure}

To begin, pick $x \in K \cap oq_i$ such that $d(x,o) \gg M$.  Then
\begin{equation} \label{equation (*)2}
d(q_i, x) - d(q_i, o) =  - d (x, o).
\end{equation}
Approximate the four point metric space $\{ x,o, q_i, p_i \}$ by a tree.  Let the corresponding four points in the
tree be $\{ \tilde{x}, \tilde{o}, \tilde{q_i}, \tilde{p_i} \}$.  We may assume this tree is at worst $(1,
2\delta)$-quasi-isometric to $\{x,o, q_i, p_i \}$ \cite[pg.33]{GH} (see Figure \ref{b}).  We thus obtain the
inequalities
\begin{eqnarray*}
d(p_i, x) - d(p_i, o) & \ge & d(\tilde{p_i}, \tilde{x}) - d(\tilde{p_i}, \tilde{o}) - 4\delta  \\ & = &
d(\tilde{x}, \tilde{o}) - 2 ( \tilde{x} \, | \, \tilde{p_i}) - 4\delta \\ &=& d(\tilde{x}, \tilde{o}) - 2 (
\tilde{q_i} \, | \, \tilde{p_i}) - 4\delta \\ & \ge & d(x,o) - 2\delta - 2 (q_i \, | \, p_i ) - 6 \delta - 4\delta
\\
    & \ge & d(x,o) - 2M - 12 \delta > \frac{d(x,o)}{2}.
\end{eqnarray*}
Together with equation (\ref{equation (*)2}), this contradicts equation (\ref{equation (*)}).
\end{pf}

\begin{lem}
$\pi$ is continuous and surjective.
\end{lem}
\begin{pf}
To see that $\pi$ is surjective, pick a geodesic ray $\gamma$ based at $o$.  By compactness there is a sequence
$\{ t_i \} \subset (0,\infty)$ such that $\gamma (t_i)$ converges to a point $\xi \in \mathcal{H}X$.  Then by
definition $\pi (\xi) = [\gamma]$.

%As a corollary of the above lemmas, we know that if $p_n \longrightarrow \xi \in \mathcal{H}X$ then $p_n \longrightarrow \pi (\xi) \in \partial X$.  From this it follows that $\pi$ is surjective.

Let $\xi^n \longrightarrow \xi$ in $\mathcal{H}X$.  Pick sequences $\{ p_i^n \}$ such that for all $n$, $\lim_{i
\rightarrow \infty} p_i^n = \xi^n$.  This implies $\lim_{i \rightarrow \infty} p^n_i = \pi (\xi^n)$ in $\partial
X$ for all $n$.  There is an increasing sequence $i_1, i_2, i_3, \ldots$ of natural numbers such that $$\lim_{n
\rightarrow \infty} p^n_{i_n} = \xi \text{  in } \mathcal{H}X, \quad \text{and} \quad
        ( \pi (\xi^n) \, | \, p^n_{i_n}) > n.$$
This implies $\lim_{n \rightarrow \infty} p^n_{i_n} = \pi(\xi)$ in $\partial X$.  By the definition of
$\delta$-hyperbolicity $$( \pi(\xi) \, | \, \pi (\xi^n) ) \ge \left[ \text{min} \{( \pi(\xi) \, | \, p^n_{i_n} ),
\, (\pi(\xi^n) \, | \,p^n_{i_n}) \}
     - 3\delta \right] \longrightarrow \infty \quad \text{as} \quad n \rightarrow \infty.$$
Therefore $\pi (\xi^n) \longrightarrow \pi (\xi)$ in $\partial X$.
\end{pf}

The proof of the following lemma is trivial and has been omitted.

\begin{lem}
For all $g \in \text{Isom}(X)$, $g \circ \pi = \pi \circ g$.
\end{lem}

\begin{lem}
\label{atomless} If $X$ admits a cocompact isometric action, then $\pi_* \mu_p$ has no atoms for all $p \in X$.
\end{lem}
\begin{pf}
As $\pi_* \mu_p \ll \pi_* \mu_o$, it is sufficient to show $\pi_* \mu_o$ has no atoms.  Suppose there exists
$\alpha \in \partial X$ such that $\pi_* \mu_o (\alpha) > 0$.  Since $X$ admits a cocompact isometric action,
there is a constant $D>0$ such that any $p \in X$ is at most a distance $D$ from an orbit point $g_p.o$ for $g_p
\in \text{Isom}(X)$.  It follows that the total mass of the measure $\mu_p$ is at most $e^{D h(X)}$.

We first show all the horofunctions in the fiber $\pi^{-1} (\alpha)$ are a bounded distance from each other.  Pick
a geodesic ray in $X$ based at $o$ asymptotic to $\alpha$.  Pick a sequence $p_i$ of points going to infinity on
the geodesic ray.  After passing to a subsequence we may assume the points $p_i$ converge to some $\xi \in
\pi^{-1} (\alpha) \subset \mathcal{H}X$.

By definition $$b_\xi (x) = \lim \left( -2 (p_i \, | \, x) + d(x,o) \right).$$ We know $$(\alpha \, | \, x ) -
2\delta \le \liminf (p_i \, | \, x ) \le
    (\alpha \, | \, x ).$$
Because $b_\xi$ is well defined, the $\liminf$ can be replaced by a limit.  Thus $$-2 (\alpha \, | \, x ) + d(x,o)
\le b_\xi (x) \le -2(\alpha \, | \, x ) +
    d(x, o ) + 4\delta.$$
So for any other $\zeta \in \pi^{-1} (\alpha)$ we have $$| b_\xi (x) - b_\zeta (x) | \le 4 \delta.$$

Using this we obtain the inequality
\begin{eqnarray*}
\pi_* \mu_{p_i} (\alpha) = \int_{\pi^{-1} (\alpha) }
    d\mu_{p_i} &=& \int_{\pi^{-1} (\alpha)} e^{-h(X) b_\zeta
        (p_i)} d\mu_o (\zeta) \\
&\ge & \int_{\pi^{-1} (\alpha)} e^{-h(X) [ b_\xi (p_i)
        + 4 \delta ]} d\mu_o (\zeta) \\
& = &  e^{-h(X)[ b_\xi (p_i)+4 \delta]} \cdot \pi_* \mu_o (\alpha)>0 .
\end{eqnarray*}
But $\lim_{i \rightarrow \infty} b_\xi (p_i) = -\infty.$ So $$\lim_{i \rightarrow \infty} \pi_* \mu_{p_i} (\alpha)
= \infty.$$ Since $\pi_* \mu_{p_i}$ has total mass at most $e^{D h(X)}$, this is a contradiction.  Therefore
$\pi_* \mu_o$ has no atoms.
\end{pf}

\section{the Besson-Courtois-Gallot inequality}
\label{inequality}

\begin{thm}
\label{mainthm} Let ${M_{\text{hyp}}}$ be a closed hyperbolic $n$-manifold for $n \ge 3$.  Let $Z$ be a compact
$n$-dimensional almost everywhere Riemannian metric space with universal cover $X$. (For a definition, see Section
\ref{aers}.)  Let \mbox{$f: Z \longrightarrow M_\h$} be a homotopy equivalence.  Then $$h(X)^n \, \text{Vol}(Z)
\ge (n-1)^n \text{Vol} ({M_{\text{hyp}}}).$$
\end{thm}

Let $f$ also denote the lifted map $f: X \longrightarrow \Hn$.  Recall that if $Z$ and $M_\h$ are compact, then a
homotopy equivalence $f: Z \longrightarrow M_\h$ lifts to a quasi-isometry $f: X \longrightarrow \U{M_\h}$ between
the universal covers.  This implies that the volume growth entropy of $X$ is a strictly positive and finite (see
for example \cite[Prop.5.10]{Gr}).  Fix a basepoint $o \in X$, and let $f(o) \in \Hn$ be a basepoint of $\Hn$.
Since $f: X \longrightarrow \Hn$ is a quasi-isometry, $f$ extends to a homeomorphism between the Gromov boundaries
$f: \partial X \longrightarrow \partial \Hn$.

We now define the natural map $F: X \longrightarrow \Hn$, but postpone proving its regularity properties until
later sections.  For $x \in X$, let $\mu_x$ be the Patterson-Sullivan measure at $x$.  Recall we have defined a
continuous map $\pi: \mathcal{H}X \longrightarrow \partial X$.  Push forward the Patterson-Sullivan measure
$\mu_x$ on $\mathcal{H}X$ to a probability measure $(f \circ \pi)_* \mu_x$ on $\partial \Hn$.  Define $F: X
\longrightarrow \Hn$ by $$F(x) := \text{barycenter} ((f \circ \pi)_* \mu_x).$$ $F$ is the \emph{natural map
induced by $f$}.  $F$ is a $\Gamma := \pi_1 (Z) \cong \pi_1 (M_\h)$-equivariant continuous map.  It therefore
descends to a continuous map $F: Z \longrightarrow M_{\text{hyp}}$.

\begin{rmk} \label{symmetries}
Any additional symmetries of the map $f: Z \longrightarrow M_\h$ also become symmetries of $F$.  Namely, if $Z$
and $M_\h$ possess an isometric involution, and $f: Z \longrightarrow M_\h$ is equivariant with respect to the
involutions, then $F: Z \longrightarrow M_\h$ is similarly equivariant.  This follows immediately from the
definition, because both $\pi: \mathcal{H}X \longrightarrow \partial X$ and Patterson-Sullivan measure are
Isom$(X)$-equivariant.  This fact will be used in the proofs of Theorems \ref{cores thm} and \ref{n-mfds}.
\end{rmk}

$F,f: X \longrightarrow \Hn$ are $\Gamma$-equivariant maps.  The straight-line homotopy between them is also
$\Gamma$-equivariant.  Therefore the downstairs maps $F,f : Z \longrightarrow M_\h$ are homotopic.  This shows $F:
Z \longrightarrow M_\h$ is a homotopy equivalence.  As $M_\h$ is a closed manifold, $F$ must be surjective.

Proving the local regularity properties of the natural map requires some work.  These properties are summarized in

\begin{prop}
\label{mainprop} The natural map $F: Z \longrightarrow {M_{\text{hyp}}}$ has the following properties:\\ (1)  It
is locally Lipschitz and differentiable almost everywhere.\\ (2)  $|\text{Jac} F (p) | \le \left( \frac{h(X)}{n-1}
\right)^n$
    almost everywhere.\\
(3)  If for some $p$, $|\text{Jac} F (p) | = \left( \frac{h(X)}{n-1} \right)^n$, then the differential $dF_p$ is a
homothety of ratio $\left( \frac{h(X)}{n-1} \right)^n$.
\end{prop}
{\noindent}This proposition will be proven in later sections.  Specifically, (1) will be proven in Section
\ref{lipschitz}, (2) and (3) will be proven in Section \ref{jacobian}.  Let us now temporarily assume it, and
complete the proof of Theorem \ref{mainthm}.  (In this section, we use only (1) and (2).  (3) will not be used
until Section \ref{equality}.)

By assumption, $Z$ is an almost everywhere Riemannian metric space.  So by definition $Z$ has a subset $\Omega$ of
full measure admitting a $\mathcal{C}^1$-atlas $\{ U_\phi, \phi \}_{\phi \in \Phi}$ and a
$\mathcal{C}^0$-Riemannian metric $\{ g_\phi \}_{\phi \in \Phi}$.  This Riemannian metric induces a volume element
$\omega_Z$ which agrees with $n$-dimensional Hausdorff measure on $Z$.  Let $\omega_{M_\h}$ be the volume element
on ${M_\h}$.

\begin{lem} \label{dagger}
\begin{eqnarray*}
\text{Vol}(M_\h) = \text{Vol}(F(Z)) =
        \int_{F(Z)} \omega_{M_\h}
    &\le & \int_Z |\text{Jac} F | \ \omega_Z \\
&\le & \left( \frac{h(X)}{n-1} \right)^n \text{Vol}(Z). \end{eqnarray*}
\end{lem}

\begin{pf}
Assuming Proposition \ref{mainprop}, the only non-trivial part is to prove $$\int_{F(Z)} \omega_{M_\h}
    \le  \int_Z |\text{Jac} F | \ \omega_Z.$$
This amounts to justifying the change of variables formula for the singular space $Z$.  To do this, we will unpack
the definitions and apply the change of variables formula for Lipschitz maps.

$Z \setminus \Omega$ is measure zero, and $F$ is locally Lipschitz.  This implies $F(Z \setminus \Omega) \subset
M_\h$ is also measure zero.  The collection of open sets $\{ U_\phi \}_{\phi \in \Phi}$ covers $\Omega$.  Let $\{
E_k \}$ be a countable partition of $\cup_\phi U_\phi \subset Z$ into measurable sets such that each $E_k$ is
contained in some open set $U_\phi$.  Since $F(\cup_k E_k) = F( \cup_\phi U_\phi) \subset M_\h$ is of full
measure, it suffices to prove the inequality on each measureable set $E_k$, i.e. it suffices to show that
$$\int_{F(E_k)} \omega_{M_\h} \le  \int_{E_k} |\text{Jac} F | \ \omega_Z.$$ Assume that $E_k \subset U_\phi$.  Let
us also assume without a loss of generality that the image $F(U_\phi) \subset M_\h$ lies in an open set $V_M$
equipped with a smooth diffeomorphism $\psi$ onto an open subset of $\mathbb{R}^n$.  Let $g_M$ denote the smooth
Riemannian metric on $\psi(V_M)$ given by the hyperbolic metric on $M_\h$.

The volume element $(\phi^{-1})^* ({\omega_Z}|_{U_\phi})$ is defined on an open subset of $\mathbb{R}^n$.  It is
induced by the the $\mathcal{C}^0$-Riemannian metric $g_\phi$.  Concretely this means that $$(\phi^{-1})^*
({\omega_Z}|_{U_\phi}) = \sqrt{\det{g_\phi}} \, dx_1 dx_2 \ldots dx_n.$$ Similarly, $$(\psi^{-1})^*
{\omega_{M_\h}} = \sqrt{\det{g_M}} \, dy_1 dy_2 \ldots dy_n.$$ Define the locally Lipschitz map $G:=\psi \circ F
\circ \phi^{-1}$, which is a map between subsets of Euclidean space.  By definition, the Jacobian of $F: Z
\longrightarrow M_\h$ at $p \in U_\phi$ is $$|\text{Jac}F(p)| := \frac{\sqrt{\det{g_M} (\psi
(F(p)))}}{\sqrt{\det{g_\phi} (\phi(p))}}
        \cdot |\text{Jac}G (\phi(p))|.$$
Applying the change of variables formula for Lipschitz maps \cite[3.4.3]{EG} to $G$ yields \begin{eqnarray*}
\int_{\phi(U_\phi)} & \chi_{\phi(E_k)}(x) & \sqrt{ \det{g_M} (G(x))}   |\text{Jac}G(x)| \ dx_1 dx_2 \ldots dx_n \\
 & & = \int_{\psi(V_M)} \#{\{G^{-1}(y) \cap \phi(E_k) \}} \ \sqrt{ \det{g_M} (y)} \ dy_1 dy_2 \ldots dy_n .
\end{eqnarray*}
Since $\det{g_\phi}$ vanishes on a set of measure zero, way may perform the following step
\begin{eqnarray*}
 \int_{\phi(U_\phi)} && \chi_{\phi(E_k)}(x) \  \sqrt{ \det{g_M} (G(x))} \ |\text{Jac}G(x)| \ dx_1 dx_2 \ldots dx_n  \\
 & =& \int_{\phi(U_\phi)} \chi_{\phi(E_k)}(x) \ \frac{\sqrt{\det{g_M} (G(x))}}{\sqrt{\det{g_\phi} (x)}} \
        \sqrt{\det{g_\phi}(x)} \   |\text{Jac}G(x)| \ dx_1 dx_2 \ldots dx_n  \\
 & =& \int_{\phi(U_\phi)} \chi_{\phi(E_k)}(x) \ |\text{Jac}F(\phi^{-1}(x))| \ \sqrt{\det{g_\phi}(x)} \  dx_1 dx_2 \ldots dx_n.
\end{eqnarray*}
Putting this together yields
\begin{eqnarray*}
\int_{E_k} |\text{Jac}F(z)| \ \omega_Z & = & \int_{\phi(U_\phi)} \chi_{\phi(E_k)}(x) \ |\text{Jac}F(\phi^{-1}(x))|
\ \sqrt{\det{g_\phi}(x)} \  dx_1 dx_2 \ldots dx_n \\ &=&  \int_{\psi(V_M)} \#{\{G^{-1}(y) \cap \phi(E_k) \}} \
\sqrt{ \det{g_M} (y)} \ dy_1 dy_2 \ldots dy_n  \\ &= & \int_{V_M} \# {\{F^{-1}(m) \cap E_k \}} \ \omega_{M_\h}(m)
\ge
   \int_{F(E_k)} \omega_{M_\h}.
\end{eqnarray*}
  \end{pf}

The above lemma implies that $$h(X)^n \, \text{Vol}(Z)  \ge  (n-1)^n \, \text{Vol}({M_\h}). $$ This proves the
inequality of Theorem \ref{mainthm}.  \square

\section{The Barycenter map is locally Lipschitz}
\label{lipschitz}

This section proves part (1) of Proposition \ref{mainprop}.

\vskip 3pt {\noindent}(1)  The natural map $F: Z \longrightarrow M_\h$ is locally Lipschitz and differentiable
almost everywhere. \vskip 3pt {\noindent}It is only necessary to prove $F$ is locally Lipschitz.  Almost
everywhere differentiability will then follow by using Rademacher's theorem (see Section \ref{aers}).  We will
prove the lifted map $F: X \longrightarrow \Hn$ is locally Lipschitz by factoring it as a composition of two
locally Lipschitz maps.  Namely we will define a locally Lipschitz map $\Phi : X \longrightarrow L^2 (\mathcal{H}
X)$ and a $\mathcal{C}^1$-map $P : L^2 (\mathcal{H} X) \longrightarrow \Hn$ such that $F = P \circ \Phi$.

The barycenter map takes a positive atomless measure $\nu$ on $\partial \Hn$ to the unique point $x = \text{bar}
(\nu)$ defined implicitly by the equation $$\int_{\partial \Hn} \langle \nabla B^o_\theta, \, v \rangle_x \, d\nu
(\theta) = 0 \quad \text{for all } v\in T_x \Hn,$$ where $B^o$ is the Busemann function on $\Hn$ (normalized so
$B^o(o,\theta)=0$ for all $\theta \in \partial \Hn$).

Consider the Hilbert space $L^2 (\mathcal{H}X)$ of square integrable functions on $\mathcal{H}X$ with respect to
the Patterson-Sullivan probability measure $\mu_o$.  Define a $\Gamma := \pi_1 (Z)$-action on $L^2 (\mathcal{H}X)$
by $$(\gamma. \phi)(\eta) := \phi(\gamma^{-1}.\eta) \cdot \sqrt{\exp{(-h(X) b_\eta (\gamma.o))}}.$$

\begin{lem}
$\Gamma$ acts by isometries on $L^2 (\mathcal{H}X)$.
\end{lem}
\begin{pf} \begin{eqnarray*}
\int_{\mathcal{H}X} (\gamma. \phi)^2 (\eta) \ d\mu_o (\eta) &=&
    \int_{\mathcal{H}X} \phi^2 (\gamma^{-1}.\eta) e^{-h(X) b_\eta(\gamma.o)} \ d\mu_o(\eta) \\
&=& \int_{\mathcal{H}X} \phi^2 (\gamma^{-1}.\eta) \ d\mu_{\gamma.o} (\eta) \\ &=&     \int_{\mathcal{H}X} \phi^2
(\eta) \ d(\gamma^{-1}_* \mu_{\gamma .o} )(\eta)
    = \int_{\mathcal{H}X} \phi^2 \ d\mu_o
\end{eqnarray*}
\end{pf}

Let $L^2_+ (\mathcal{H}X)$ denote the strictly positive functions in $L^2 (\mathcal{H}X)$.  Notice that $\Gamma$
acts by isometries on $L^2_+ (\mathcal{H}X)$.  An element $\phi \in L^2_+ (\mathcal{H}X)$ defines a positive
atomless measure $\phi^2 \, d\mu_o$.  Push this measure forward via the map $f \circ \pi$ to a measure $(f \circ
\pi)_* (\phi^2 d\mu_o)$ on $\partial \Hn$.  Define the map $P: L^2_+ (\mathcal{H}X) \longrightarrow \Hn$ by $$P :
\phi \longmapsto \text{bar} ((f \circ \pi)_* (\phi^2 d\mu_o)).$$ In other words  $P(\phi)$ is the unique point $x$
defined implicitly by the equation $$\int_{\partial \Hn} \langle \nabla B^o_\theta, \, v \rangle_x \, d((f \circ
\pi)_* (\phi^2 d\mu_o))(\theta)
     = \int_{\mathcal{H} X} \langle \nabla B^o_{ f \circ \pi (\eta)}, \, v \rangle_x \, \phi^2 (\eta) d\mu_o = 0,$$
for all $v \in T_x \Hn$.

\begin{lem}
$P$ is $\Gamma$-equivariant.
\end{lem}
\begin{pf}
$P(\gamma.\phi)$ is the unique point $x$ such that for all $v \in T_x \Hn$,
\begin{eqnarray*}
0 &=& \int_{\mathcal{H}X}
    \langle \nabla B^o_{ f \circ \pi (\eta)}, \, v \rangle_x \
  \phi^2 (\gamma^{-1} \eta) e^{-h(X) b_\eta (\gamma.o)} d\mu_o (\eta) \\
&=& \int_{\mathcal{H}X}
    \langle \nabla B^o_{f \circ \pi (\eta)}, \, v \rangle_x \
    \phi^2 (\gamma^{-1} \eta) d\mu_{\gamma.o} (\eta) \\
&=& \int_{\mathcal{H}X}
   \langle \nabla B^o_{ f \circ \pi \circ \gamma (\eta)}, \, v \rangle_x \
   \phi^2 (\eta) d\mu_o(\eta) \qquad \text{  (change of variables)} \\
&=& \int_{\mathcal{H}X}
   \langle \nabla B^o_{ \gamma \circ f \circ \pi (\eta)}, \, v \rangle_x \
\phi^2 (\eta)d\mu_o (\eta)\qquad\text{  ($f \circ \pi$ is $\Gamma$-equiv.)}\\ &=& \int_{\mathcal{H}X} \langle
\nabla B^o_{f \circ \pi (\eta)}, \,
    d \gamma^{-1} (v) \rangle_{\gamma^{-1}.x} \
    \phi^2 (\eta) d\mu_o (\eta) \quad \text{ ($\nabla B^o$ is Isom$(\Hn)$ inv.)}
\end{eqnarray*}
Since $d \gamma^{-1} : T_x \Hn \longrightarrow T_{\gamma^{-1}.x} \Hn$ is an isomorphism, this implies
$\gamma^{-1}.x = P (\phi)$.  Therefore $P(\gamma. \phi) = x=\gamma. P(\phi)$.
\end{pf}

Pick a $\mathcal{C}^\infty$ frame $\{ e_i \}$ on $\Hn$ and define a map \mbox{$Q: \Hn \times L^2_+ (\mathcal{H}X)
\longrightarrow \mathbb{R}^n$} by
\begin{eqnarray*}
\lefteqn{Q: (x, \phi) \longmapsto} & & \\ & &   \left( \int_{\mathcal{H}X} \langle \nabla B^o_{f \circ \pi
(\eta)}, \, e_1 \rangle_x \,
      \phi^2 (\eta) d\mu_o (\eta), \ldots, \int_{\mathcal{H}X}
      \langle \nabla B^o_{f \circ \pi (\eta)}, \, e_n \rangle_x \, \phi^2 (\eta) d\mu_o (\eta)
           \right).
\end{eqnarray*}
$\mathcal{H}X$ is compact, $B^o_\theta (x)$ is a $\mathcal{C}^\infty$ function of both $\theta$ and $x$, and $\{
e_i \}$ is a $\mathcal{C}^\infty$ frame.  Using these facts, applying the Lebesgue dominated convergence theorem
proves that $Q$ is $\mathcal{C}^\infty$.  Notice $P$ is defined implicitly by the equation $$Q (P(\phi), \phi) =
(0, \ldots, 0).$$ The goal is to show that $P$ is $\mathcal{C}^1 $.  This can be accomplished by employing the
implicit function theorem.  (The implicit function theorem is true on Banach spaces.  See \cite[pg.366]{RS}.)  For
each fixed $\phi$, we must show the map $Q^\phi : x \longmapsto Q(x,\phi)$ has an invertible differential at each
point $x$ of the fiber $(Q^\phi)^{-1} (0, \ldots, 0)$.  Split $Q^\phi$ into coordinate functions
$Q^\phi=(Q^\phi_1, \ldots, Q^\phi_n)$.  Then
\begin{eqnarray*}
\frac{\partial}{\partial x_j} Q^\phi_i (x) &=&
   \int_{\mathcal{H}X}  (\Hess B^o_{f \circ \pi (\eta)})_x (e_j, e_i)
         \, \phi^2(\eta) \  d\mu_o (\eta) \\
& & + \int_{\mathcal{H}X} \langle \nabla B^o_{f \circ \pi (\eta)},
    \nabla_{e_j} e_i \rangle_x \, \phi^2 (\eta) \ d\mu_o (\eta).
\end{eqnarray*}
The second term in this sum satisfies
\begin{eqnarray*}
\lefteqn{\int_{\mathcal{H}X} \langle \nabla B^o_{f \circ \pi (\eta)},
    \nabla_{e_j} e_i \rangle_x \, \phi^2 (\eta) \ d\mu_o (\eta)} \\
&=& \int_{\mathcal{H}X} \langle \nabla B^o_{f \circ \pi (\eta)} \, , \,
    \sum_{k=1}^n c_{ji}^k e_k \rangle_x \, \phi^2 (\eta)\ d\mu_o(\eta)
= \sum_{k=1}^n c_{ji}^k Q^\phi_k (x) = 0,
\end{eqnarray*}
for some constants $c_{ji}^k$ depending on the frame $\{ e_i \}$.  This implies the bilinear form on $T_x \Hn$
determined by the differential of $Q^\phi$ at $x$ satisfies
\begin{eqnarray*}
\langle v, dQ^\phi_x (v) \rangle  &=&
    \int_{\mathcal{H}X} (\Hess B^o_{f \circ \pi (\eta)})_x (v, v)
        \ \phi^2(\eta)  d\mu_o (\eta)  \\
    &= & \int_{\partial \Hn} (\Hess B^o_\theta)_x (v, v) \
        d \left((f \circ \pi)_* \, (\phi^2 d\mu_o)\right) (\theta),
\end{eqnarray*}
for all $v \in T_x \Hn$.  The right hand side of this equation is strictly positive by Lemma \ref{positive
definite}.  This implies the differential of $Q^\phi$ at $x$ is positive definite, and thus invertible.  Therefore
the implicit function theorem may be applied to conclude that $P$ is $\mathcal{C}^1$.

Define a map
\begin{eqnarray*}
\Phi: X  &\longrightarrow  & L^2_+ (\mathcal{H}X) \subset L^2 (\mathcal{H}X) \\
  x & \longmapsto & \sqrt{\exp{(- h(X) b_\eta (x))}}.
\end{eqnarray*}

\begin{lem}
$\Phi$ is Lipschitz.
\end{lem}
\begin{pf}
Let $D$ be the diameter of the downstairs metric space $Z$ covered by $X$.  It follows that for any $p \in X$, $\|
\Phi(p) \|$ is at most $e^{D h(X)}$.  Pick points $x,y \in X$.  The goal is to control the quantity
$$\int_{\mathcal{H}X} |e^{-\frac{1}{2} h(X) b_\eta (x)}
    - e^{-\frac{1}{2} h(X) b_\eta (y)} |^2 d\mu_o (\eta).$$
For $\eta \in \mathcal{H} X$, $b_\eta$ is $1$-Lipschitz.  This implies $ b_\eta (y)  \le  d(y,x) + b_\eta (x).$
Using this we obtain the inequalities
\begin{eqnarray*}
e^{-\frac{1}{2} h(X) b_\eta (x)} - e^{-\frac{1}{2} h(X) b_\eta (y)} & \le &
    e^{-\frac{1}{2} h(X) b_\eta (x)} - e^{-\frac{1}{2} h(X) d(y,x)}
    e^{-\frac{1}{2} h(X) b_\eta (x)} \\
&=&     e^{-\frac{1}{2} h(X) b_\eta (x)} ( 1 - e^{-\frac{1}{2} h(X) d(y,x)}). \\ \Rightarrow |e^{-\frac{1}{2} h(X)
b_\eta (x)} -e^{-\frac{1}{2} h(X) b_\eta (y)}| &\le& (e^{-\frac{1}{2} h(X) b_\eta (x)}+e^{-\frac{1}{2} h(X) b_\eta
(y)})
    (1- e^{-\frac{1}{2} h(X) d(x,y)}).
\end{eqnarray*}
Therefore,
\begin{eqnarray*}
\lefteqn{\int_{\mathcal{H}X} | e^{-\frac{1}{2} h(X) b_\eta (x)} -
    e^{-\frac{1}{2} h(X) b_\eta (y)}|^2 d\mu_o (\eta)} & & \\
&\le & (1- e^{-\frac{1}{2} h(X) d(x,y)})^2 \int_{\mathcal{H}X}
        |e^{-\frac{1}{2} h(X) b_\eta (x)} +e^{-\frac{1}{2} h(X) b_\eta (y)}|^2 d\mu_o (\eta) \\
&= & (1- e^{-\frac{1}{2} h(X) d(x,y)})^2 \cdot \| \Phi(x) + \Phi(y) \|^2 \\ &\le & (1- e^{-\frac{1}{2} h(X)
d(x,y)})^2 \cdot \left[
    \| \Phi(x) \| + \| \Phi(y) \| \right]^2
 \le  (1- e^{-\frac{1}{2} h(X) d(x,y)})^2 \cdot 4 e^{2D h(X)}
\end{eqnarray*}
To complete the proof, notice that for $t \ge 0$, \mbox{$(1-e^{-\frac{1}{2} h(X) t}) \le \frac{1}{2} h(X)  t$}.
Applying this yields $$\int_{\mathcal{H}X} |e^{-\frac{1}{2} h(X) b(x, \eta)}
    - e^{-\frac{1}{2} h(X) b_\eta (y)} |^2 d\mu_o (\eta) \le
    \left( \frac{1}{2} h(X) d(x,y) \right)^2 \cdot 4 e^{2D h(X)}.$$
So finally we've obtained $$\frac{\| \Phi(x) - \Phi(y) \|}{d(x,y)} \le h(X) \cdot e^{D h(X)}.$$
\end{pf}

$P$ is $\mathcal{C}^1$, therefore $F = P \circ \Phi$ is locally Lipschitz.  $\Gamma$-equivariance implies $F$
descends to a locally Lipschitz map $X/\Gamma = Z \longrightarrow \Hn/\Gamma = M_\h$.  By Rademacher's theorem,
$F$ is differentiable almost everywhere (see Section \ref{aers}).

\section{The Jacobian Estimate}
\label{jacobian}

This section proves parts (2) and (3) of Proposition \ref{mainprop}.

\vskip 3pt {\noindent}(2)  $|\text{Jac} F (p) | \le \left( \frac{h(X)}{n-1} \right)^n$
    almost everywhere.\\
(3)  If for some $p$, $|\text{Jac} F (p) | = \left( \frac{h(X)}{n-1} \right)^n$, then the differential $dF_p$ is a
homothety of ratio $\left( \frac{h(X)}{n-1} \right)^n$. \\ \vskip 3pt

{\noindent}The proof closely follows Section 5 of \cite{F}, which is in turn based on
\cite[pgs.636-639]{BCGergodic}.  Recall that $X$ is the universal cover of $Z$, $\Gamma := \pi_1 (Z) \cong \pi_1
(M_\h)$, $\{e_i \}$ is a $\mathcal{C}^\infty$ frame on $T \Hn$, and $\Omega \subset X$ is a subset of full measure
possessing a $\mathcal{C}^1$-atlas (see Sections \ref{gdrs} and \ref{aers}).

In Section \ref{lipschitz}, we defined a $\Gamma$-equivariant locally Lipschitz map $\Phi: X \longrightarrow L^2_+
(\mathcal{H}X)$, a $\Gamma$-equivariant $\mathcal{C}^1$ barycenter map $P:L^2_+ (\mathcal{H}X) \longrightarrow
\Hn$, and a $\mathcal{C}^\infty$ map $Q = (Q_1, \ldots, Q_n) : \Hn \times L^2_+ (\mathcal{H}X) \longrightarrow
\mathbb{R}^n$.  They satisfied the equations $F = P \circ \Phi$ and $Q(P(\phi), \phi) = 0$ for $\phi \in L^2_+
(\mathcal{H}X)$.  We thus obtain the implicit equation $Q (F , \Phi ) = 0.$ Let $\mathcal{O} \subseteq \Omega
\subseteq X$ be the set of points where $\Phi$ is differentiable.  By Rademacher's theorem (see Section
\ref{aers}), $\mathcal{O} \subseteq X$ is a subset of full measure. Pick $p \in \mathcal{O}$ and $v \in T_p X$.
\begin{lem}\label{derivative}
Differentiating the function $$x \longmapsto Q_i (F(x), \Phi(x)) = \int_{\mathcal{H}X} \langle \nabla B^o_{f \circ
\pi (\eta)}, \, e_i
        \rangle_{F(x)} \, e^{- h(X) b_\eta (x)} \, d\mu_o (\eta) \ = \, 0$$
at $p$ in the direction of $v$ yields
\begin{eqnarray*}
0 &=& \int_{\mathcal{H}X} (\Hess B^o_{f \circ \pi (\eta)} )_{F(p)}
    (dF (v), e_i ) \  e^{-h(X) b_\eta (p)} d\mu_o (\eta) \\
  & & + \int_{\mathcal{H}X} \langle \nabla B^o_{f \circ \pi(\eta)}
    , \, \nabla_{dF (v)} e_i \rangle_{F(p)}
        \  e^{-h(X) b_\eta (p)} d\mu_o (\eta) \\
  & & + \int_{\mathcal{H} X} \langle \nabla B^o_{f \circ \pi(\eta)}, \, e_i \rangle_{F(p)} \
   [-h(X) \, \langle \nabla b_\eta, \, v \rangle_p \, e^{-h(X) b_\eta (p)}] \ d \mu_o (\eta). \qquad (\star)
\end{eqnarray*}
\end{lem}
\begin{pf}
This is an application of the Lebesgue dominated convergence theorem.  $\mathcal{H}X$ is compact.  So to apply the
theorem it is sufficient to find a $c>0$ such that for all $\eta \in \mathcal{H}X$, the function $$x \longmapsto
\langle \nabla B^o_{f \circ \pi (\eta)} , \, e_i \rangle_{F(x)} \,  e^{-h(X) b_\eta (x)} $$ is locally
$c$-Lipschitz near $p$.  To show the existence of such a constant, use that $b_\eta$ is $1$-Lipschitz,
$\mathcal{H}X$ is compact, and $\langle \nabla B^o_\theta, \, e_i \rangle_q : \Hn \times \partial \Hn
\longrightarrow \mathbb{R}$ is $\mathcal{C}^\infty$.
\end{pf}
{\noindent}The second term in equation $(\star)$ satisfies
\begin{eqnarray*}
\lefteqn{\int_{\mathcal{H}X} \langle \nabla B^o_{f \circ \pi(\eta)} , \, \nabla_{dF (v)} e_i \rangle_{F(p)} \
e^{-h(X) b_\eta (p)} d\mu_o (\eta)  } \\ &=&  \int_{\mathcal{H}X} \langle \nabla B^o_{f \circ \pi(\eta)}, \,
\sum_{j=1}^n c_{ij}  e_j \rangle_{F(p)} \  e^{-h(X) b_\eta (p)} d\mu_o (\eta) \\ &=& \sum_{j=1}^n c_{ij} \, Q_j
(F(p), \Phi(p)) = 0,
\end{eqnarray*}
for some constants $c_{ij}$ depending on the frame $\{ e_i \}$.  Therefore, for $v \in T_p X$, $u \in T_{F(p)}
\Hn$ we have
\begin{eqnarray*}
\lefteqn{\int_{\mathcal{H}X} (\Hess B^o_{f \circ \pi (\eta)} )_{F(p)}
    (dF (v), u) \  e^{-h(X) b_\eta (p)} d\mu_o (\eta) } \\
&=&  h(X) \int_{\mathcal{H} X} \langle \nabla B^o_{f \circ \pi(\eta)}, \, u \rangle_{F(p)} \ \langle \nabla
b_\eta, \, v \rangle_p \, e^{-h(X) b_\eta (p)} \, d \mu_o (\eta). \qquad (\star \star)
\end{eqnarray*}

Let $\| \mu_p \|$ denote the total mass $\mu_p (\mathcal{H} X)$ of the measure $\mu_p$.  Define a bilinear form
$K$ and a quadratic form $H$ on $T_{F(p)} \Hn$ by $$\langle K (w) \, , \, u \rangle_{F(p)} := \frac{1}{\| \mu_p
\|} \,
  \int_{\mathcal{H}X } (\Hess B^o_{f \circ \pi(\eta)} )_{F(p)} (w,u)
    \, e^{-h(X) b_\eta (p)} d \mu_o (\eta) $$
$$ \langle H(u) \, , \, u \rangle_{F(p)} := \frac{1}{\| \mu_p \|} \,
 \int_{\mathcal{H}X} \langle \nabla B^o_{f \circ \pi (\eta)}, \, u \rangle^2_{F(p)} \,
     e^{-h(X) b_\eta (p)} d \mu_o (\eta).$$
Note that the symmetric endomorphism $K$ is positive definite by Theorem \ref{positive definite}.  It is therefore
invertible.  This is used in the proof of Lemma \ref{bounding the jacobian}.

Use equation $(\star \star)$, the Cauchy-Schwarz inequality, and the definition of $K$ to obtain
\begin{eqnarray*}
\lefteqn{| \langle K \circ dF (v), u \rangle_{F(p)} | } \\ & & = \left| \frac{h(X)}{\| \mu_p \|} \int_{\mathcal{H}
X} \langle \nabla B^o_{f \circ \pi(\eta)} , \, u \rangle_{F(p)}
  \cdot \langle \nabla b_\eta , \, v \rangle_p \cdot e^{-h(X) b_\eta (p)} \  d \mu_o (\eta) \right| \\
& & \le \frac{h(X)}{\| \mu_p \|}  \left[ \int_{\mathcal{H} X}
    \langle \nabla B^o_{f \circ \pi (\eta)}, \, u \rangle_{F(p)}^2 \, d\mu_p (\eta)\right]^{1/2}
   \cdot \left[ \int_{\mathcal{H} X}
      \langle \nabla b_\eta, \, v \rangle_p^2 \ d\mu_p (\eta) \right]^{1/2} \\
& & = \frac{h(X)}{\| \mu_p \|^{1/2}} \ \left[ \langle H(u), u \rangle_{F(p)} \right]^{1/2} \cdot
   \left[ \int_{\mathcal{H} X}
      \langle \nabla b_\eta, \, v \rangle_p^2 \ d\mu_p (\eta) \right]^{1/2},
\end{eqnarray*}
for all $u \in T_{F(p)} \Hn$ and $v \in T_p \Omega$.

\begin{lem}\label{bounding the jacobian} \cite[pg.637]{BCGergodic}
For all $p \in \mathcal{O}$, $$|\text{Jac} F(p) | \le \frac{h(X)^n (\det H)^{1/2}}{n^{n/2} \det K} .$$
\end{lem}
\begin{pf}
The proof is the proof of Lemma 5.4 of \cite[pg.637]{BCGergodic} with two modifications.  First replace the
Busemann function $B_\alpha$ with the horofunction $b_\eta$, and notice that $$\sum_{i=1}^n \langle \nabla b_\eta,
v_i \rangle_p^2 = \|\nabla \, b_\eta \|^2_p \le 1. \qquad (\ddagger)$$ Second, the total mass of the
Patterson-Sullivan measure $\| \mu_p \|$ must be carried through the estimate, but it cancels itself out in the
final step.
\end{pf}

A key property of hyperbolic space is that Busemann functions on $\Hn$ satisfy the equation $$(\Hess
B^o_\theta)(u,v) = \langle u,v \rangle - \langle \nabla B^o_\theta, \, u \rangle \cdot \langle \nabla B^o_\theta,
\, v \rangle$$ for all $\theta \in \partial \Hn$ (see \cite[pgs.750-751]{BCGlong}).  Integrating this equation
over $\mathcal{H} X$ yields $$K = \text{Id} - H.$$

\begin{lem}
The symmetric endomorphism $H$ is positive definite.
\end{lem}
\begin{pf}
Suppose there exists $x \in X$ and $u \in T_x X$ such that $\langle H(u), u \rangle_x = 0$.  From the definition
of $H$, this implies the support of the measure $(f \circ \pi)_* \mu_0$ on $\partial \mathbb{H}^n$ is contained in
a codimension one conformally round sphere in $S \subset \mathbb{H}^n$.  By the equivariance of the
Patterson-Sullivan measures, the support of $(f \circ \pi)_* \mu_{\gamma.o}$ is contained in $\gamma(S)$ for
$\gamma \in \Gamma$.  But $\mu_{\gamma.o} \ll \mu_o$ implies the support of $(f \circ \pi)_* \mu_{\gamma.o}$ is
contained in $S$.  Therefore $\Gamma$ preserves $S \subset \partial \mathbb{H}^n$.  This contradicts the fact that
the limit set of $\Gamma$ is the entire sphere at infinity.
\end{pf}

A short computation shows that $\text{tr}(H) =1$.  Therefore we can apply the brain in a jar lemma (Lemma
\ref{brain in a jar lemma}) to $H$.  This yields $$\frac{\det(H)}{[\det (\text{Id} - H)]^2} \le
    \left[ \frac{n}{(n-1)^2} \right]^n,$$
with equality if and only if $H = \frac{1}{n} \text{Id}$.  Combining this with Lemma \ref{bounding the jacobian}
proves the desired inequality, namely $$|\text{Jac} F(p)| \le \left( \frac{h(X)}{n-1} \right)^n.$$ This completes
the proof of part (2) of Proposition \ref{mainprop}.

The proof of part (3) of Proposition \ref{mainprop} is the ``equality case'' argument on page 639 of
\cite{BCGergodic} with two modifications (as in Lemma \ref{bounding the jacobian}).  First replace the Busemann
function $B_\alpha$ with the horofunction $b_\eta$ and use inequality $(\ddagger)$.  Second, the total mass of the
Patterson-Sullivan $\| \mu_p \|$ must be carried through the estimate, and again it cancels itself out.  This
completes the proof of Proposition \ref{mainprop}.

\section{The equality case}
\label{equality}

\begin{thm}
\label{mainthm2} Let ${M_{\text{hyp}}}$ be a closed hyperbolic $n$-manifold for $n \ge 3$.  Let $Z$ be a compact
$n$-dimensional convex Riemannian amalgam (see Section \ref{amalgams}).  Let \mbox{$f: Z \longrightarrow M_\h$} be
a homotopy equivalence.  If $$h(X)^n \, \text{Vol}(Z) = (n-1)^n \text{Vol} ({M_{\text{hyp}}})$$ then the natural
map $F: Z \longrightarrow M_\h$ is a homothetic homeomorphism.
\end{thm}

\begin{rmk}
Recall that if $Z$ is either \\ $\bullet$  the metric doubling of a convex hyperbolic manifold with boundary (see
Section \ref{crmwb}), or \\ $\bullet$  a cone-manifold (see Section \ref{cone-manifolds}), \\ then $Z$ is a convex
Riemannian amalgam.
\end{rmk}

Let $X$ be the universal cover of $Z$.  Lift the convex Riemannian amalgam structure on $Z$ to a convex Riemannian
amalgam structure on $X$ with decomposition $\{ C_j \}$ into convex Riemannian manifolds with boundary.  Define
the incomplete disconnected Riemannian manifold $\Omega := \bigcup_j \text{int}(C_j) \subseteq X$.

The equation $$h(X)^n \, \text{Vol}(Z) = (n-1)^n \, \text{Vol} ({M_{\text{hyp}}})$$ implies the string of
inequalities from Lemma \ref{dagger} is in fact a string of equalities.  Therefore $|\text{Jac} F (p) | = \left(
\frac{h(X)}{n-1} \right)^n$ almost everywhere.  By Proposition \ref{mainprop}, $dF$ is a homothety almost
everywhere.  The goal is to prove $F: Z \longrightarrow M_\h$ is a homothetic homeomorphism.  Without a loss of
generality, scale the metric of $Z$ so that $dF$ is an isometry (not necessarily orientation preserving) almost
everywhere.  This forces $\text{Vol}(Z) = \text{Vol}({M_{\text{hyp}}})$.  It now suffices to show that $F: Z
\longrightarrow M_\h$ is an isometric homeomorphism.  This will be done by working on the universal covers and
showing $F: X \longrightarrow \Hn$ is an equivariant isometric homeomorphism.

\begin{lem}
\label{contraction} $F: X \longrightarrow \Hn$ is a contraction mapping, i.e. for any pair $x,y \in X$, $d_{\Hn}
(F(x),F(y)) \le d_X (x,y)$.
\end{lem}
\begin{pf}
Pick a length minimizing geodesic segment $\alpha$ connecting $x$ to $y$.  Define $\mathcal{O} \subseteq \Omega
\subseteq X$ to be the set of full measure where $dF$ is an isometry.  There exists a small perturbation
$\alpha_\varepsilon$ of $\alpha$ such that $\alpha_\varepsilon \cap \mathcal{O} \subseteq \alpha_\varepsilon$ is a
subset of full measure, and length$(\alpha_\varepsilon) \le \text{length}(\alpha) + \varepsilon$.  $F$ is locally
Lipschitz, so in all Hausdorff dimensions $F$ maps measure zero sets to measure zero sets.  This implies $F$
preserves length on $\alpha_\varepsilon$.  Therefore $$d_{\Hn} (F(x), F(y)) \le
\text{length}(F(\alpha_\varepsilon)) =
    \text{length} (\alpha_\varepsilon)
    \le \text{length} (\alpha) + \varepsilon = d_X (x,y) + \varepsilon.$$
Since $\varepsilon$ was arbitrary, the result follows.
\end{pf}

We claim that $F$ is volume preserving.  $| \text{Jac} F | = 1 \text{ a.e.}$ and $F$ is locally Lipschitz, so $F$
is volume non-increasing.  If it strictly decreased the volume of some measurable set, then we would have Vol$(Z)>
\text{Vol}({M_{\text{hyp}}})$.  This would be a contradiction.  Therefore $F$ is volume preserving.

There exist constants $C_1,\varepsilon_\h > 0$ such that if $y \in \Hn$ and $\varepsilon < \varepsilon_\h$ then
\begin{eqnarray}\label{1}
v_n \varepsilon^n  \le
    \text{Vol} (B_{\Hn} (y, \varepsilon)) \le v_n \varepsilon^n
        (1 + C_1 \varepsilon^2),
\end{eqnarray}
where $v_n$ is the volume of a unit ball in $\mathbb{R}^n$. Similarly, $\Omega \subseteq X$ is a Riemannian
manifold with sectional curvature bounded from above.  (An upper curvature bound follows from the definition of a
convex Riemannian amalgam.)  For compact $K \subset \Omega$ (possibly a point), define $$\text{inj}_\Omega (K) :=
\inf_{z \in K} (\text{injectivity radius of $\Omega$ at $z$}).$$ The upper curvature bound implies there exist
constants $C_2, \varepsilon_\Omega >0$ such that if $\varepsilon < \varepsilon_\Omega$, $\varepsilon <
\text{inj}_\Omega (z)$, and $B_X (z, \varepsilon) \subset \Omega$ then
\begin{eqnarray}\label{2}
v_n \varepsilon^n (1- C_2 \varepsilon^2) \le \text{Vol}(B_\Omega (z, \varepsilon)).
\end{eqnarray}

We now consider the restriction of $F$ to the ``smooth'' set $\Omega$.  Define $F_\Omega := F|_\Omega: \Omega
\longrightarrow \Hn$.

\begin{lem}
\label{BCGlemma1} $F_\Omega$ is injective.
\end{lem}
\begin{pf}
(This proof is an adaptation of \cite[Lem.C.4]{BCGlong}.)  Suppose there exist distinct $z_1, z_2 \in \Omega$ such
that $F_\Omega (z_1) = F_\Omega (z_2)= y.$  Pick $$\varepsilon_0 < \min \{\varepsilon_\Omega,\ \text{inj}_\Omega
(z_1), \
    \text{inj}_\Omega (z_2), \ \varepsilon_\h \}$$
sufficiently small such that $$B_X (z_1, \varepsilon_0) \cap B_X (z_2, \varepsilon_0) = \emptyset \quad \text{and}
\quad
        B_X (z_1, \varepsilon_0) \cup B_X (z_2, \varepsilon_0) \subset \Omega.$$
In particular, inequality \ref{1} (resp. inequality \ref{2}) is valid at $y$ (resp. $z_1$ and $z_2$) for
$\varepsilon < \varepsilon_0$.

Since $F$ is contracting, $F_\Omega (B_\Omega (z_i, \varepsilon)) \subseteq
        B_{\Hn} (y, \varepsilon)$ for all $\varepsilon < \varepsilon_0$.  Therefore
$$B_\Omega (z_1, \varepsilon)
    \bigcup B_\Omega (z_2, \varepsilon) \subseteq F_\Omega^{-1} (B_{\Hn} (y,\varepsilon)).$$
Since $F$ is volume preserving (and $X \setminus \Omega$ is measure zero),
\begin{eqnarray*}
\text{Vol}(B_{\Hn} (y,\varepsilon)) = \text{Vol}(F_\Omega^{-1} (B_{\Hn} (y, \varepsilon))) &\ge & \text{Vol}
\left( B_\Omega (z_1, \varepsilon) \bigcup
    B_\Omega (z_2, \varepsilon)\right) \\
    &=& \text{Vol} (B_\Omega (z_1, \varepsilon)) + \text{Vol}
        (B_\Omega (z_2, \varepsilon)).
\end{eqnarray*}
Now apply inequalities \ref{1} and \ref{2} to obtain $$v_n \varepsilon^n (1 + C_1 \varepsilon^2) \ge
    2 v_n \varepsilon^n ( 1- C_2 \varepsilon^2).$$
This implies $\varepsilon^2 (C_1 + 2 C_2) \ge 1$ for all $\varepsilon < \varepsilon_0$, which is a contradiction.
\end{pf}

Since $\Omega$ is locally compact, Lemma \ref{BCGlemma1} shows that $F_\Omega: \Omega \longrightarrow
F_\Omega(\Omega)$ is a homeomorphism.

\begin{lem}
\label{BCGlemma2} $F_\Omega^{-1}$ is locally Lipschitz.
\end{lem}
\begin{pf}
(This proof is an adaptation of \cite[Lem.C.7]{BCGlong}.)  Pick $z \in \Omega$.  Let $y := F_\Omega (z)$.  Pick
$\varepsilon_0 <  \min \{\varepsilon_\h , \varepsilon_\Omega \}$ sufficiently small such that $B_X (z, 3
\varepsilon_0) \subset \Omega$, $\varepsilon_0 < \text{inj}_\Omega (\, \overline{B_X (z, 2 \varepsilon_0) }\,)$,
and $$ 2 \varepsilon_0^2 (C_1 + C_2)  < \frac{1}{2^n}.$$ Define $V:= B_\Omega (z, \varepsilon_0)$.  Since $F$ is a
contraction mapping, $F(V) \subseteq B_{\Hn} (y, \varepsilon_0)$.  We will show that if $z_1,z_2 \in V$ are
distinct, then $$d_X (z_1,z_2) < 2 d_{\Hn} (F(z_1), F(z_2)).$$ This will imply ${F|_V}^{-1}$ is $2$-Lipschitz on
the open set $F(V)$.

Suppose the above inequality is false, i.e. there exist distinct $z_1,z_2 \in V$ such that $2 d_{\Hn} (F(z_1),
F(z_2)) \le d_X (z_1,z_2).$  Set $\varepsilon := d_{\Hn} (F(z_1), F(z_2)) \le \varepsilon_0.$  Since $\varepsilon
<  \varepsilon_0$, inequality \ref{1} remains valid at $F(z_1)$ and $F(z_2)$ for $\varepsilon$.   $B_{\Hn}
(F(z_1), \varepsilon)$ and $B_{\Hn} (F(z_2), \varepsilon)$ intersect, and their intersection contains a ball of
radius $\varepsilon/2$ centered in the middle of the minimizing geodesic joining $F(z_1)$ to $F(z_2)$ (see Figure
\ref{big}).  Therefore $$\text{Vol} \left( B_{\Hn} (F(z_1), \varepsilon) \bigcup B_{\Hn} (F(z_2), \varepsilon)
    \right) \le v_n \left[ 2 \varepsilon^n (1 + C_1 \varepsilon^2) -
    \frac{1}{2^n} \varepsilon^n  \right].$$
Moreover, $d_X (z_1, z_2) \ge 2\varepsilon$ implies $B_X (z_1, \varepsilon) \cap B_X (z_2, \varepsilon) =
\emptyset$.  Since\\ \mbox{$B_X (z_1, \varepsilon) \cup B_X (z_2, \varepsilon) \subset B(z, 2 \varepsilon_0)
\subset \Omega$}, we may apply inequality \ref{2} to conclude
\begin{eqnarray*}
\text{Vol}\left( F \left( B_X (z_1, \varepsilon) \bigcup B_X (z_2, \varepsilon) \right) \right) &=& \text{Vol}
\left( B_X (z_1, \varepsilon) \bigcup B_X (z_2, \varepsilon) \right) \\ & \ge & v_n [ 2 \varepsilon^n ( 1 - C_2
\varepsilon^2) ].
\end{eqnarray*}

\begin{figure}[ht]
\begin{center}
\epsfig{file=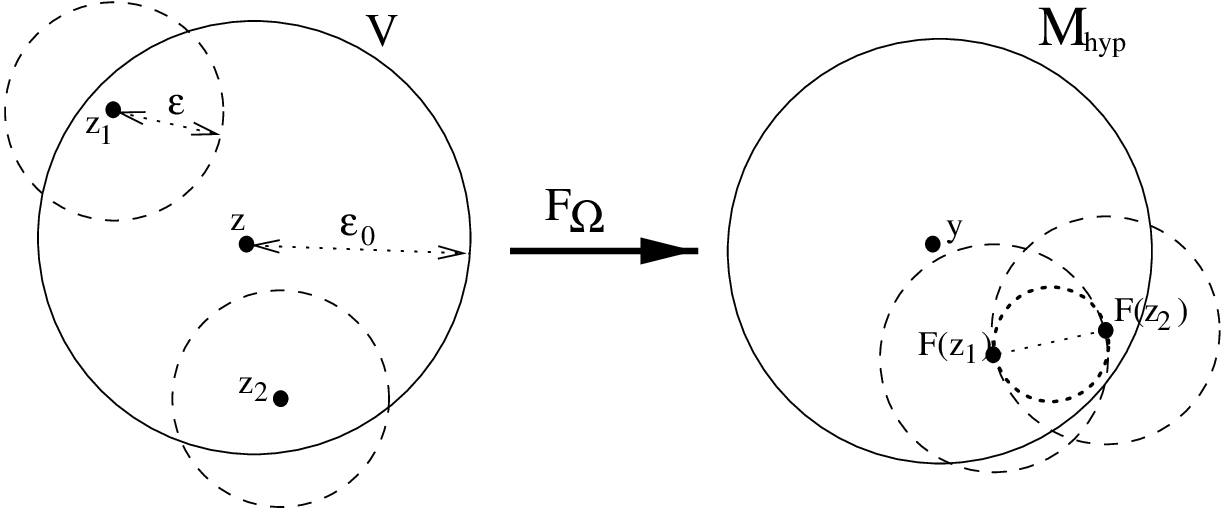,angle=-90} \caption{} \label{big}
\end{center}
\end{figure}

$F$ is a contraction mapping.  Therefore
\begin{eqnarray*}
F(B_X (z_1, \varepsilon)) \subseteq B_{\Hn} (F(z_1), \varepsilon)
    & \text{  and  } &
F(B_X (z_2, \varepsilon)) \subseteq B_{\Hn} (F(z_2), \varepsilon).
\end{eqnarray*}
Putting these inequalities together yields $$ v_n \left[ 2 \varepsilon^n (1 + C_1 \varepsilon^2) -
    \frac{1}{2^n} \varepsilon^n  \right]
    \ge v_n \varepsilon^n [2 (1- C_2 \varepsilon^2)],$$
implying $$ 2 \varepsilon^2 (C_1 + C_2)   \ge \frac{1}{2^n}.$$ This contradicts our choice of $\varepsilon_0$.
\end{pf}

Since ${F|_\Omega}^{-1}$ is locally Lipschitz, it is differentiable almost everywhere.  Therefore for almost every
$z \in \Omega$ $$d(Id_\Omega)_z = d( {F|_\Omega}^{-1} \circ F_\Omega )_z
    =d({F|_\Omega}^{-1}) \circ dF_\Omega.$$
This implies $d({F|_\Omega}^{-1})$ is an isometry of tangent spaces almost everywhere.  By an argument analogous
to Lemma \ref{contraction}, and by working on small balls in ${M_{\text{hyp}}}$, one can see that
${F|_\Omega}^{-1}$ is locally a contraction mapping $F(\Omega) \longrightarrow \Omega$.  Both $F_\Omega$ and
${F|_\Omega}^{-1}$ are local contraction maps.  Thus they are both local isometries.

\begin{lem}
For all $j$, $F|_{C_j}$ is an isometry onto its image, and $F(\text{int}(C_j))$ is convex.
\end{lem}
\begin{pf}
By continuity, it suffices to show $F|_{\text{int} (C_j)}$ is an isometry onto its image.  Pick $x,y \in
\text{int}(C_j)$.  Since $\text{int}(C_j)$ is convex, there exists a geodesic segment $xy \subset \text{int}(C_j)$
joining $x$ to $y$.  $F_\Omega$ is a local isometry.  Thus length$(F (xy)) = \text{length}(xy)$, and $F(xy)
\subset \Hn$ is locally a hyperbolic geodesic joining $F(x)$ to $F(y)$.  In $\Hn$, a local geodesic is a global
geodesic.  Therefore $$d_{\Hn} (F(x), F(y)) = \text{length} (F (xy)) = \text{length} (xy) = d_X (x,y).$$
\end{pf}

\begin{lem}
\label{balls} For each $z \in X$, there exists $\delta_z >0$ such that $d_X (z,x) < \delta_z$ implies $d_{\Hn}
(F(z), F(x)) = d_X (z,x)$.
\end{lem}
{\noindent}(Notice this lemma proves neither that $F$ is an isometry nor that it is locally injective.  To see
why, consider the branched double covering of an Euclidean disk by an Euclidean cone with cone angle $4\pi$.)

\begin{pf}
Pick $z \in X$.  If $z \in \Omega$ then we are done.  So assume $z \in X \setminus \Omega$.  Up to rearranging the
indices we may assume $$z \in \left( \partial C_1 \cup \partial C_2 \cup \ldots \cup \partial C_l \right )
\setminus \left( \partial C_{l+1} \cup \ldots \cup \partial C_N \right).$$ By this, there exists $\delta_z > 0$
such that $$B_X (z, \delta_z) \cap \Omega \subseteq \bigcup_{j=1}^l C_j.$$ For $x \in B_X (z, \delta_z)$, there
exists an integer $j_x$ and a geodesic segment $xz \subset C_{j_x}$.  $F|_{C_{j_x}}$ is an isometry, so $d_X (z,x)
= d_{C_{j_x}} (z,x) = d_{\Hn} (F(z), F(x))$.  This proves the lemma.
\end{pf}

As a corollary of this lemma, we see that $F^{-1} (y)$ is a discrete set for all $y \in \Hn$.

\begin{lem}
$F: X \longrightarrow \Hn$ is injective.
\end{lem}
\begin{pf}
For $y \in \Hn$, $F^{-1} (y)$ is the discrete set $\{ z_i \}$.  We will show $F^{-1} (y)$ must be a single point.
For all $i$, pick $\delta_i < 1$ such that $d_X (z_i,x) < \delta_i$ implies $d_{\Hn} (F(z_i), F(x)) = d_X
(z_i,x)$.  Assume the metric balls $B_X (z_i, \delta_i)$ are disjoint.

If for some $i$, $F|_{B_X(z_i, \delta_i)} : B_X(z_i, \delta_i) \longrightarrow B_{\Hn} (F(z_i), \delta_i)$ is not
surjective, then $F$ can be properly homotoped to a map taking $B_X (z_i, \delta_i) \longrightarrow \partial
B_{\Hn} (F(z_i), \delta_i)$.  Moreover, this can be done without altering $F$ outside of $B_X (z_i, \delta_i)$.

If for every $i$, $F|_{B_X(z_i, \delta_i)} : B_X(z_i, \delta_i) \longrightarrow B_{\Hn} (F(z_i), \delta_i)$ is not
surjective, then $F$ is properly homotopic to a map which does not have $y$ in its image.  $F$ is a proper
surjective homotopy equivalence.  Thus every map properly homotopic to $F$ is surjective.  Therefore for some $i$,
$F|_{B_X(z_i, \delta_i)}$ is surjective.  We may assume $i=1$.

$F_\Omega$ is injective and open, $\Omega$ is open and dense, and $F(B_X (z_1, \delta_1)) = B_{\Hn}(y, \delta_1)$.
Thus $$F(\Omega \setminus B_X (z_1, \delta_1)) \bigcap B_{\Hn} (y, \delta_1)
    = \emptyset.$$
By extending to the closure of $\Omega$, this shows that for all $i>1$, $$F(z_i) = y \notin  B_{\Hn} (y,
\delta_1).$$  Therefore $F^{-1} (y)$ must be the single point $z_1$.
\end{pf}

Recall that $F:X \longrightarrow \Hn$ is surjective.  Therefore $F: X \longrightarrow \Hn$ is a continuous
bijection.  Since $X$ is locally compact, $F$ is a homeomorphism.  To prove that $F$ is an isometry, it is
sufficient to show $F^{-1}$ is a contraction mapping.

\begin{lem}
$F^{-1}$ is a contraction mapping, i.e. it is $1$-Lipschitz.
\end{lem}
\begin{pf}
By the above lemmas, $F$ imposes a convex Riemannian amalgam structure on $\Hn$.  The collection of convex
Riemannian manifolds with boundary is $\{ F(C_j) \}$.

Pick $x,y \in \Hn$, $\varepsilon > 0$.  Since each $F(C_j)$ is convex, there exists a path $\gamma \subset \Hn$
joining $x$ to $y$ such that: \\ $\bullet $ length$(\gamma) \le d(x,y) + \varepsilon$, \\ $\bullet $ $\gamma \cap
F(\partial C_j)$ is at most two points for any $F(C_j)$ of the decomposition. \\ As the collection $\{ F(C_j) \}$
is locally finite, the set $F^{-1} (\gamma ) \setminus \Omega$ is finite.  $F|_{C_j}$ is an isometry for all $j$.
Therefore $F^{-1} (\gamma)$ is a path of the same length as $\gamma$.  This implies $F^{-1}$ is $1$-Lipschitz.
\end{pf}

This completes the proof of Theorem \ref{mainthm2}.

\section{Applications}

\subsection{Kleinian groups}\label{kleinian section}

The Kleinian group theory notation used here is defined in Sections \ref{convex cores}-\ref{deformation theory
section}.

Let $N$ be a compact acylindrical $3$-manifold (see Section \ref{pared 3-manifolds section}).  Recall that by
Corollary \ref{pared acylindrical}, there exists a convex cocompact hyperbolic $3$-manifold $M_g$ such that
$C_{M_g}$ is homeomorphic to $N$ and the boundary of the convex core $\partial C_{M_g} \subset M_g$ is totally
geodesic.

As was discussed in Section \ref{introduction}, the following theorem solves a conjecture in Kleinian groups.
\begin{thm}
\label{cores thm} Let $N$ be a compact acylindrical $3$-manifold.  Let $(M_g, m_g)$ be a convex cocompact
hyperbolic $3$-manifold such that $C_{M_g}$ is homeomorphic to $N$ and the boundary of the convex core $\partial
C_{M_g} \subset M_g$ is totally geodesic.  For all $(M,m) \in \text{H}(N)$, $$\text{Vol}(C_M) \ge \text{Vol}
(C_{M_g}),$$ with equality if and only if $M$ and $M_g$ are isometric.
\end{thm}

Fix an $(M,m) \in \text{H}(N)$.  We may assume without a loss of generality that $M$ is geometrically finite.
Since $\partial C_{M_g}$ is totally geodesic, metrically doubling $C_{M_g}$ across its boundary produces a compact
hyperbolic manifold $DC_{M_g}$.

\begin{lem}
\label{Alexandrov lemma} Metrically doubling the convex core $C_M$ across its boundary yields an Alexandrov space
with curvature bounded below by $-1$.
\end{lem}
\begin{pf}
In \cite[Appendix A]{S} it was proven that taking an $\varepsilon$-neighborhood of $C_M$ in $M$, and metrically
doubling that across its boundary to obtain $D\mathcal{N}_\varepsilon C_M$, yields an Alexandrov space with
curvature bounded below by $-1$.  $DC_M$ is the Gromov-Hausdorff limit of these spaces as $\varepsilon \rightarrow
0$.  Being an Alexandrov space with curvature bounded below by $-1$ is a closed property in the Gromov-Hausdorff
topology \cite[Prop.10.7.1,pg.376]{BBI}.  This proves the lemma.
\end{pf}

\begin{rmk}
Lemma \ref{Alexandrov lemma} also follows from a more general theorem of Perelman (see Theorem \ref{doubling
theorem}).
\end{rmk}

{\noindent}In particular, by Lemma \ref{Alexandrov lemma} and Theorem \ref{Perel'man}, the volume growth entropy
of $\U{DC_M}$ is not greater than $2 = h(\mathbb{H}^3)$.  This will be used later.  Let $\sigma$ denote the
boundary preserving isometric involution of $DC_M$ and $DC_{M_g}$.

\vskip 6pt {\noindent}\textbf{Case 1:} Assume $M$ is convex cocompact.\\

\textit{Proof of Case 1:} Both $M$ nor $M_g$ are convex cocompact.  By Theorem \ref{homeo thm}, $m_g \circ m^{-1}$
is homotopic to a homeomorphism $C_M  \longrightarrow C_{M_g}$.  ($m^{-1}$ denotes a relative homotopy inverse of
$m$.)  This homeomorphism can be ``doubled'' to produce a homeomorphism between the doubled manifolds $f: DC_M
\longrightarrow DC_{M_g}$.  Theorem \ref{mainthm} may now be applied to $f: {DC_M} \longrightarrow {DC_{M_g}}$.
This proves that $$h(\U{DC_M})^3 \ \text{Vol}(DC_M) \ge 2^3 \ \text{Vol}(DC_{M_g}),$$ with equality if and only if
$DC_M$ and $DC_{M_g}$ are isometric. Since $h(\U{DC_M}) \le 2$ we have $$\text{Vol}(DC_M) \ge
\text{Vol}(DC_{M_g}),$$ with equality if and only if $DC_M$ and $DC_{M_g}$ are isometric. Dividing both sides by
$2$ yields the desired inequality.

Let us now assume $\text{Vol}(DC_M) = \text{Vol}(DC_{M_g})$.  The goal is to show $M$ and $M_g$ are isometric.  To
do this it is sufficient to show $C_M$ and $C_{M_g}$ are isometric.  The map $f: DC_M \longrightarrow DC_{M_g}$ is
by construction $\sigma$-equivariant.  Let $F: DC_M \longrightarrow DC_{M_g}$ be the natural map induced by $f$.
The $\sigma$-equivariance of $f$ implies $F$ is also $\sigma$-equivariant (see Remark \ref{symmetries}).
Therefore $F: C_M \longrightarrow C_{M_g}$ is an isometry.  This completes the proof of \textit{Case 1}. \square

\vskip 6pt {\noindent}\textbf{Case 2:} $(M,m)$ is not convex cocompact.\\

\textit{Proof of Case 2:}  $(M,m) \in \text{H}(N)$, so by definition $m: N \longrightarrow M^o$ is a homotopy
equivalence (see Section \ref{deformation theory section}).  Moreover, by Theorem \ref{homeo thm} there exists a
homeomorphism $g: N \longrightarrow C_M \cap M^o$.  Let $P := g^{-1} (C_M \cap \partial M^o) \subset \partial N$.
Consider $(M,m)$ as an element of $\text{H}(N,P)$ with no additional parabolics.  $(N,P)$ is pared acylindrical.
Therefore there exists a geometrically finite $M_P \in \text{H}(N,P)$ with no additional parabolics such that
$C_{M_P}$ has totally geodesic boundary.

Again by Theorem \ref{homeo thm}, there is a homeomorphism $C_{M_P} \longrightarrow N \setminus P$, inducing a
doubled homeomorphism $DC_{M_P} \longrightarrow D(N\setminus P)$ between open manifolds.  Let $p \subset P \subset
\partial N$ be a finite collection of disjoint simple closed curves such that $p$ is a strong deformation retract
of $P$, i.e. each component of $p$ is a core curve of a component of $P$.  Then there is a homeomorphism
$D(N\setminus P) \longrightarrow (DN) \setminus p$. Moreover, the manifolds $DC_{M_P} \cong D(N \setminus P)$ are
topologically obtained by removing the curves $p \subset \partial N \subset DN$. Conversely, replacing the removed
curves of  $DN \setminus p$ corresponds to performing a (topological) Dehn surgery on $DN \setminus p$. Therefore
the homeomorphism type of $DN \cong DC_{M_g}$ can be obtained by performing a topological Dehn surgery on
$DC_{M_P}$. $DC_{M_P}$ is a finite volume hyperbolic manifold.  Therefore Mostow rigidity implies that $DC_{M_g}$
may in fact be obtained by performing a hyperbolic Dehn surgery on $DC_{M_P}$ \cite{Th}.

Hyperbolic Dehn surgery strictly decreases volume \cite[Thm.6.5.6]{Th}\cite{Bes}.  Therefore $\text{Vol}(DC_{M_P})
> \text{Vol}(DC_{M_g}).$  Moreover, by \cite{NZ} there exists a closed hyperbolic $3$-manifold $L$ obtained from
hyperbolic Dehn surgery on $DC_{M_P}$ such that $$\text{Vol}(DC_{M_P}) > \text{Vol}(L) > \text{Vol}(DC_{M_g}).$$
So to complete the proof it suffices to show that $$\text{Vol} (DC_M) \ge \text{Vol}(L).$$  This will be
accomplished by geometrically filling in the cusps and reducing to the case of closed manifolds.

The geometric (though not hyperbolic) Dehn surgery arguments below are based on techniques from \cite{Bes,L}. (See
also \cite{BCS1} for another application of these methods.)  The exposition here will roughly follow \cite{Bes}.

Outside a compact set, the metric on $DC_M$ is a collection of smooth rank two hyperbolic cusps.  For simplicity,
let us assume that $DC_M$ has exactly one cusp.  The general case follows by performing the following geometric
operations on each cusp individually.

Pick a compact exhaustion $K_i$ of $DC_M$ such that each boundary $\partial K_i$ is a smooth horospherical torus.
By ``opening up'' the cusp of $DC_M$ one can construct a family of metric spaces $\{ (A,d_{\varepsilon})
\}_{\varepsilon \in (0, 1]}$ such that: \\
 1.  $A$ is homeomorphic to $K_i$. \\
 2. $(A,d_\varepsilon)$ is an Alexandrov space with curvature bounded below by $-1 - c_\varepsilon^2$, where $\lim_{\varepsilon \rightarrow 0}
c_\varepsilon = 0$. \\
 3.  For each $i$, there is an isometric embedding $K_i \longrightarrow (A, d_\varepsilon)$
for all sufficiently small $\varepsilon$. \\
 4.  $\lim_{\varepsilon \rightarrow 0} \text{Vol}(A, d_\varepsilon) = \text{Vol}(DC_M)$. \\
 5.  Near the boundary of $A$ the metric $d_\varepsilon$ is a flat Riemannian metric with totally geodesic boundary $\partial A$. \\
 6.  For any $\varepsilon \le 1$, $(\partial A, d_\varepsilon) = (\partial A, \varepsilon \cdot d_1)$.

{\noindent}(A careful and clear explanation of this procedure is in \cite[Sec.2.2]{Bes}.  See also \cite{L}.)

Let $W$ denote a solid torus.  Using the fact that $W$ is a product of a disk and a circle, one can build a
product Riemannian metric $g$ on $W$ such that: \\ 1.  $(W, g)$ has totally geodesic boundary isometric to an
Euclidean torus. \\ 2.  $(W, g)$ has sectional curvature bounded below by zero.

A manifold homemorphic to $L$ is obtained from $A$ by appropriately gluing $\partial W$ to $\partial A$ (i.e. by
an appropriate Dehn surgery).  The goal is to perform this gluing geometrically.  To do so we must scale
appropriately and interpolate between the flat torus boundaries of $(A, d_\varepsilon)$ and $(W,g)$.  Let $T^2
\times [0,1]$ denote a trivial interval bundle over a torus.  Pick diffeomorphisms $$\phi_0: T^2 \times \{0 \}
\longrightarrow \partial A \text{  and  }
   \phi_1 : T^2 \times \{ 1 \} \longrightarrow \partial W$$
such that the glued up manifold $$A \ \bigcup_{\phi_0} \ (T^2 \times [0,1]) \ \bigcup_{\phi_1}\ W$$ is
homeomorphic to $L$.

Consider the metrics $\phi_0^* d_\varepsilon$ and $\phi_1^* g$ on the boundary of $T^2 \times [0,1]$.  We now
apply a lemma from \cite{Bes}.

\begin{lem} \cite[App.A.2]{Bes}
For any $n>0$ there exist $\alpha_n, \varepsilon_n >0$ and a Riemannian metric $\sigma_n$ on $T^2 \times [0,1]$
such that:\\
 1.  The curvature of $\sigma_n$ is bounded between $-1/n$ and $1/n$. \\
 2.  The volume of $(T^2 \times [0,1], \sigma_n)$ is less than $1/n$. \\
 3.  $(T^2 \times [0,1], \sigma_n)$ has totally geodesic boundary given by $(T^2 \times \{ 0 \}, \phi_0^*
d_{\varepsilon_n})$ and $(T^2 \times \{ 1 \}, \alpha_n \cdot \phi_1^* g)$.\\
 4.  $\alpha_n$ and $\varepsilon_n$ go to zero as $n \rightarrow \infty$.
\end{lem}

We may now glue geometrically to form $$(Y, \Delta_n) := (A, d_{\varepsilon_n})\ \bigcup_{\phi_0} \ (T^2 \times
[0,1], \sigma_n) \ \bigcup_{\phi_1} \ (W, \alpha_n \cdot g), \quad \text{such that}:$$ 1.  $\text{Vol}(Y,
\Delta_n) \longrightarrow \text{Vol}(DC_M).$\\ 2.  $(Y, \Delta_n)$ is an Alexandrov space with curvature bounded
below by $-1 - c_{\varepsilon_n}^2$. \\ 3.  $Y$ is homeomorphic to the closed hyperbolic manifold $L$.

There is a sequence $\eta_n \rightarrow 1$ such that the homothetically scaled spaces $(Y, \eta_n \cdot \Delta_n)$
are Alexandrov spaces with curvature bounded below by $-1$.  Theorem \ref{mainthm} applied to $(Y, \eta_n \cdot
\Delta_n)$ and $L$ yields the inequality $$\text{Vol}(Y, \eta_n \cdot \Delta_n) \ge \text{Vol}(L).$$ Taking $n
\rightarrow \infty$ yields the desired inequality $$\text{Vol}(DC_M) \ge \text{Vol}(L).$$

\square

Using the above geometric surgery arguments, we now prove that the inequality of Theorem \ref{cores thm} holds
also for pared acylindrical manifolds.

\begin{cor} \label{pared acyl cor} Let $(N,P)$ be a compact pared acylindrical $3$-manifold.  Let $(M_g, m_g)$ be
a hyperbolic $3$-manifold such that $C_{M_g}$ is homeomorphic to $N \setminus P$ and the boundary of
the convex core $\partial C_{M_g} \subset M_g$ is totally geodesic.  For all $(M,m) \in \text{H}(N,P)$,
$$\text{Vol}(C_M) \ge \text{Vol} (C_{M_g}).$$ \end{cor} \begin{pf} $DC_{M_g}$ is a finite volume hyperbolic
manifold.  By repeating the arguments from the beginning of Case 2 above, it follows that $DC_{M_g}$ is obtained
topologically by performing a (possibly empty) set of Dehn surgeries on the manifold $DC_M$.  (Note that since
$DC_{M_g}$ is not compact, Dehn surgery is not performed on every end of $DC_M$.  The ends corresponding to the
cusps of $DC_{M_g}$ are not changed.)

By performing an infinite sequence of \emph{hyperbolic} Dehn surgeries on $DC_{M_g}$ we obtain a sequence of
\emph{closed} hyperbolic manifolds $L_k$ such that $\text{Vol}(L_k) \nearrow \text{Vol}(DC_{M_g})$ \cite{NZ}.  For
each $k$, a manifold homeomorphic to $L_k$ can be obtained from $DC_M$ by performing an appropriate topological
Dehn surgery on each end of $DC_M$.  By repeating the geometric surgery arguments of Case 2 above, the closed
manifold $Y_k$ obtained by these Dehn surgeries on $DC_M$ can be given a sequence of metrics $\delta^n_k$ such
that:\\ 1.  $\lim_{n \rightarrow \infty} \text{Vol}(Y_k, \delta^n_k) \longrightarrow \text{Vol}(DC_M).$\\ 2.
$(Y_k, \delta^n_k)$ is an Alexandrov space with curvature bounded below by $-1$. \\ 3.  $Y_k$ is homeomorphic to
the closed hyperbolic manifold $L_k$.

{\noindent}Applying Theorem \ref{mainthm} to the sequence $\{ (Y_k, \delta^n_k) \}_n$ and $L_k$ yields
$$\text{Vol}(DC_M) \ge \text{Vol}(L_k).$$ Taking $k \rightarrow \infty$ yields $$\text{Vol}(DC_M) \ge
\text{Vol}(DC_{M_g}).$$
\end{pf}

\subsection{Cone-Manifolds}
\label{cone-manifolds section}

\begin{thm}
\label{cone-manifolds theorem} Let $Z$ be compact $n$-dimensional $(n \ge 3)$ cone-manifold built with simplices
of constant curvature $K \ge -1$.  Assume all its cone angles are $\le 2\pi$.  Let $M_\h$ be a closed hyperbolic
$n$-manifold.  If $f: Z \longrightarrow M_\h$ is a homotopy equivalence then $$\text{Vol}(Z) \ge
\text{Vol}(M_\h)$$ with equality if and only if $f$ is homotopic to an isometry.
\end{thm}

\begin{pf}
Since $K \ge -1$ and all the cone angles of $Z$ are $\le 2 \pi$, this implies $Z$ is an Alexandrov space with
curvature bounded below by $-1$ \cite[pg.7]{BGP}.  Therefore, by Theorem \ref{Perel'man}, the volume growth
entropy of $\U{Z}$ is less than or equal to $(n-1)$.  Applying Theorems \ref{mainthm} and \ref{mainthm2} proves
the theorem.
\end{pf}

\subsection{Alexandrov Spaces}
\label{alexandrov spaces section}

\begin{thm}
\label{Alexandrov theorem} Let $Z$ be a compact $n$-dimensional $(n \ge 3)$ Alexandrov space with curvature
bounded below by $-1$.   Let $M_\h$ be a closed hyperbolic $n$-manifold.  If $f: Z \longrightarrow M_\h$ is a
homotopy equivalence then $$\text{Vol}(Z) \ge \text{Vol} (M_\h).$$
\end{thm}

\begin{pf}
Otsu and Shioya proved that a finite dimensional Alexandrov space with curvature bounded below by $-1$ is almost
everywhere Riemannian \cite{OS}.  By Theorem \ref{Perel'man}, the volume growth entropy of $\U{Z}$ is $\le (n-1)$.
The theorem now follows from Theorem \ref{mainthm}.
\end{pf}

\subsection{$n$-manifolds with boundary} \label{section n-mfds}

The argument used in Section \ref{kleinian section} generalizes immediately to prove a version of the
Besson-Courtois-Gallot theorem for convex Riemannian $n$-manifolds with boundary.  (For the definition of convex
Riemannian manifolds with boundary, see Section \ref{crmwb}.)  A good example of a convex Riemannian manifold with
boundary is a convex core with positive volume.

Perelman's Doubling theorem \cite[Thm.5.2]{P} is used in the proof of Theorem \ref{n-mfds}.  Unfortunately, the
well known preprint \cite{P} remains unpublished.
\begin{thm} \label{doubling theorem} \cite[Thm.5.2]{P}
Metrically doubling any Alexandrov space with curvature bounded below by $k$ across its boundary produces an
Alexandrov space with curvature bounded below by $k$.
\end{thm}

\begin{thm}
\label{n-mfds} Let $Z$ be a compact convex Riemannian $n$-manifold with boundary $(n \ge 3)$.  Assume the
sectional curvature of int$(Z)$ is bounded below by $-1$.  Let $Y_\g$ be a compact convex hyperbolic $n$-manifold
with totally geodesic boundary.   Let $f: (Z , \partial Z) \longrightarrow (Y_\g, \partial Y_\g )$ be a homotopy
equivalence of pairs.  Then $$\text{Vol} (Z) \ge \text{Vol} (Y_\g),$$ with equality if and only if $f$ is
homotopic to an isometry.
\end{thm}

\begin{pf}
The homotopy equivalence of a pairs $f: (Z , \partial Z) \longrightarrow (Y_\g, \partial Y_\g )$ can be extended
to a homotopy equivalence between the doubled manifolds $f: DZ \longrightarrow DY_\g$.  We know the sectional
curvature of int$(Z)$ is bounded below by $-1$.  So $Z$ is an Alexandrov space with curvature bounded below by
$-1$.  By Perelman's Theorem \ref{doubling theorem}, $DZ$ is also an Alexandrov space with curvature bounded below
by $-1$.  Theorem \ref{Perel'man} then implies the volume growth entropy of $\U{DZ}$ is not greater than $h(\Hn) =
(n-1)$.  Applying Theorem \ref{mainthm} yields the desired inequality.

Assume the inequality is an equality.  Theorem \ref{mainthm2} then implies the natural map $F: DZ \longrightarrow
DY_\g$ is an isometry.  As before, let $\sigma$ be the boundary preserving isometric involution of $DZ$ and
$DY_\g$.  Since $f$ is $\sigma$-equivariant, $F$ is $\sigma$-equivariant (see Remark \ref{symmetries}).
Therefore, by Theorem \ref{mainthm2}, $F: Z \longrightarrow Y_\g$ is an isometry.
\end{pf}

\begin{sc}
\noindent Stanford University \\
Mathematics, Bldg. 380 \\
450 Serra Mall \\
Stanford, CA 94305-2125 
\end{sc}

\begin{thebibliography}{BCG2}

\bibitem[AMR]{AMR}
I. Aitchison, S. Matsumotoi, J.H. Rubinstein
\textit{Immersed surfaces in cubed manifolds},
Asian J. Math. \textbf{1}, (1997), no.1, 85-95

\bibitem[Bes]{Bes}
L. Bessi{\`e}res
\textit{Sur le volume minimal des vari{\'e}t{\'e}s ouvertes},
Ann. Inst. Fourier 50(3) (2000), 965-980

\bibitem[BBI]{BBI}
D.~Burago, Y.~Burago, and S.~Ivanov
\textit{A Course in Metric Geometry},
Graduate Studies in Mathematics, Volume 33.  AMS, Providence, 2001

%\bibitem[BCG1]{BCGold}
%G.~Besson, G.~Courtois, and S.~Gallot
%\textit{Volume et entropie minimale des espaces localement sym\'etriques},
%Invent. math. \textbf{103} (1991), 417-445

\bibitem[BCG1]{BCGlong}
G.~Besson, G.~Courtois, and S.~Gallot
\textit{Entropies et rigidit\'es des espaces localement sym\'etriques de courbure strictement n\'egative},
Geom. Funct. Anal. \textbf{5}, no. 5, 731-799

\bibitem[BCG2]{BCGergodic}
G.~Besson, G.~Courtois, and S.~Gallot
\textit{Minimal entropy and Mostow's rigidity theorems},
Ergod. Th. \& Dynam. Sys. \textbf{16} (1996), 623-649

%\bibitem[BCS1]{BCS}
%J.~Boland, C. Connell, and J. Souto
%\textit{Volume rigidity for finite volume manifolds},
%preprint, (2002)

\bibitem[BCS]{BCS1}
J.~Boland, C.~Connell, and J.~Souto
\textit{Minimal volume and minimal entropy},
in \textit{Geometric structures on 3-manifolds and their deformations},
J.~Souto Clement,
Bonner Mathematische Schriften, Nr.342.  Bonn, 2001

\bibitem[BGP]{BGP}
Y.~Burago, M.~ Gromov, and G.~Perel'man
\textit{A. D. Aleksandrov spaces with curvatures bounded below},
(Russian) Uspekhi Mat. Nauk 47 (1992), no. 2(284), 3--51, 222
translation in Russian Math. Surveys 47 (1992), no. 2, 1--58

\bibitem[BGS]{BGS}
W.~Ballmann, M.~Gromov, and V.~Schroeder
\textit{Manifolds of nonpositive curvature},
Progress in Mathematics, 61.  Birkh\"auser Boston, Inc., Boston, MA, 1985

\bibitem[BM]{BM}
M. Burger and S. Mozes 
\textit{CAT(-1)-spaces, divergence groups and their commensurators},
J. Amer. Math. Soc. 9 (1996), no. 1, 57--93

\bibitem[Bon]{Bon}
F.~Bonahon
\textit{A Schl\"{a}fli-type formula for convex cores of hyperbolic 3-manifolds},
J. Differential Geom. \textbf{50} (1998), 25-58

\bibitem[BP]{BP}
R. Benedetti and C. Petronio
\textit{Lectures on Hyperbolic Geometry},
Springer-Verlag, Universitext (1992)

%\bibitem[CMT]{CMT}
%R.~Canary, Y.~Minsky, and E.~Taylor
%\textit{Spectral theory, Hausdorff dimension and the topology of hyperbolic 3-manifolds},
%Journal of Geometric Analysis \textbf{9} (1999), 17-40

%\bibitem[CEG]{CEG}
%R. Canary, D.B.A. Epstein, and P. Green
%\textit{Notes on notes of Thurston},
%in \textit{Analytical and Geometrical Aspects of Hyperbolic Spaces},
%Cambridge University Press, Cambridge (1987), 3-92

\bibitem[CHK]{CHK}
D.~Cooper, C.~Hodgson, and S.~Kerckhoff
\textit{Three-dimensional Orbifolds and Cone-Manifolds},
MSJ Memoirs, volume 5 (2000),  Mathematical Society of Japan

\bibitem[EM]{EM}
D.~Epstein and A.~Marden
\textit{Convex hulls in hyperbolic space, a theorem of Sullivan, and measured pleated surfaces}
in \textit{Analytical and Geometrical Aspects of Hyperbolic Spaces}, 
Cambridge University Press (1987), 3-92

\bibitem[EG]{EG}
L. Evans and R. Gariepy
\textit{Measure Theory and Fine Properties of Functions},
Studies in Advanced Mathematics, CRC Press (1992)

\bibitem[F]{F}
R.~Feres
\textit{The minimal entropy theorem and Mostow rigidity, after G.~Besson, G.~Courtois, and S.~Gallot},
unpublished.

\bibitem[GH]{GH}
E.~Ghys and P. de la Harpe (editors)
\textit{Sur les Groupes Hyperboliques d'apr\`es Mikhael Gromov},
Progress in Mathematics, Volume 83 (1990)

\bibitem[Gr]{Gr}
M. Gromov
\textit{Metric Structures for Riemannian and Non-Riemannian Spaces},
Progress in Mathematics, Volume 152, Birkh{"a}user, (1999)

%\bibitem[GP]{GP}
%K.~Grove and P.~Petersen (editors)
%\textit{Comparison Geometry},
%MSRI Publications, volume 30 (1997)

\bibitem[J]{J}
W.~Jaco
\textit{Lectures on Three-Manifold Topology},
CBMS Regional Conference Series in Mathematics, 43. American
Mathematical Society, Providence, R.I., (1980)

\bibitem[L]{L}
B. Leeb
\textit{$3$-manifolds with(out) metrics of nonpositive curvature},
Invent. Math. 122(2) (1995), 277-289

\bibitem[M]{M}
J.~Morgan
\textit{On Thurston's uniformization theorem for three-dimensional manifolds},
in \textit{The Smith Conjecture}, ed. by J.~Morgan and H.~Bass, Academic Press (1984), 37-125

\bibitem[NZ]{NZ}
W. Neumann and D. Zagier
\textit{Volumes of hyperbolic 3-manifolds},
Topology \textbf{24} (1985), no.3, 307-332

\bibitem[OS]{OS}
Y. Otsu and T. Shioya
\textit{The Riemannian structure of Alexandrov spaces},
J. Diff. Geom. \textbf{39} (1994), 629-658

\bibitem[P]{P}
G. Perel'man
\textit{Alexandrov spaces with curvatures bounded from below II},
preprint, (1991)

\bibitem[RS]{RS}
M. Reed and B. Simon
\textit{Functional Analysis},
Academic Press, (1980)

%\bibitem[Sp]{Sp}
%E. Spanier
%\textit{Algebraic Topology},
%Springer-Verlag, (1966)

\bibitem[S]{S}
P.~Storm
\textit{Minimal volume Alexandrov spaces},
preprint, (2001)

%\bibitem[Su]{Su}
%D. Sullivan
%\textit{Entropy, Hausdorff measures old and new, and limit sets of geometrically finite Kleinian groups},
%Acta Math. 153 (3-4) (1984), 259-277

\bibitem[Th1]{Th}
W.~Thurston
\textit{The topology and geometry of 3-manifolds},
Princeton Univ., Lecture Notes, 1976-79

\bibitem[Th2]{Th2}
W.~Thurston
\textit{Hyperbolic geometry and 3-manifolds},
in \textit{Low-dimensional topology} (Bangor, 1979), 
London Math. Soc. Lecture Note Ser., 48, 
Cambridge Univ. Press, Cambridge-New York (1982), 9-25

\bibitem[Th3]{Th3}
W.~Thurston,
\textit{Hyperbolic structures on $3$-manifolds. I. Deformation of acylindrical manifolds},
Ann. of Math. (2) 124 (1986), no. 2, 203--246

\end{thebibliography}
\end{document}